\documentclass[a4paper,12pt]{article}

\usepackage[latin1]{inputenc}

\usepackage[T1]{fontenc}
\usepackage{latexsym}
\usepackage{exscale}

\usepackage{amsmath,amsfonts,amssymb,amsthm}

\usepackage{graphicx,color}

\newtheorem{thm}{Theorem}[section]

\newtheorem{lem}[thm]{Lemma}

\newtheorem{cor}[thm]{Corollary}

\newtheorem{prop}[thm]{Proposition}

\newtheorem{rem}[thm]{Remark}

\newtheorem{dfn}[thm]{Definition}

\newcommand{\aver}[1]{-\hskip-0.46cm\int_{#1}}
\newcommand{\avert}[1]{-\hskip-0.38cm\int_{#1}}
\newcommand{\ind}{1\hspace{-2.5 mm}{1}}

\DeclareMathOperator{\supp}{supp}
\DeclareMathOperator{\Lip}{Lip}
\makeatletter
\renewcommand\@biblabel[1]{#1.}
\makeatother
\begin{document}
\allowdisplaybreaks

\author{\begin{tabular}{ccc}
Nadine Badr\footnote{Institut Camille Jordan, Universit\'e Claude Bernard, Lyon 1, UMR du CNRS 5208,
43 boulevard du 11 novembre 1918, F-69622 Villeurbanne cedex. Email: badr@math.univ-lyon1.fr} &
 Besma Ben Ali \footnote{Universit\'e de Paris-Sud, UMR du CNRS 8628, F-91405 Orsay cedex. Email:besmath@yahoo.fr}\\ \\
Universit\'e Lyon 1 & Universit\'e Paris-Sud 
\end{tabular}}

\date{\today}
\title{ $L^{p}$ Boundedness of Riesz transform related to Schr\"odinger operators on a manifold}
\maketitle
\begin{abstract}
We establish various $L^{p}$ estimates for the Schr\"odinger operator $-\Delta+V$ on Riemannian manifolds satisfying the doubling property and a Poincar\'e inequality, where $\Delta $ is the Laplace-Beltrami operator and $V$ belongs to a reverse H\"{o}lder class. At the end of this paper we apply our result on Lie groups with polynomial growth.
\end{abstract}

\tableofcontents
\section{Introduction}
The main goal of this paper is to establish the $L^{p}$ boundedness for the Riesz
transforms  $\nabla(-\Delta+V)^{-\frac{1}{2}}$, $V^{\frac{1}{2}}(-\Delta+V)^{-\frac{1}{2}}$ and related inequalities on certain classes of Riemannian manifolds. Here, $V$ is a non-negative, locally integrable function on $M$.

 For the Euclidian case, this subject was studied by many authors under different
conditions on $V$. We mention the works of Helffer-Nourrigat \cite{helffer}, Guibourg \cite{guibourg1}, Shen \cite{shen1},
Sikora \cite{sikora1}, Ouhabaz \cite{ouhabaz} and others. 

Recently, Auscher-Ben Ali \cite{auscher3} proved $L^{p}$ maximal inequalities for these operators under less
restrictive assumptions. They assumed that $V$ belongs to some reverse H\"{o}lder class $RH_{q}$ (for a definition, see section 2). A natural step further is to extend the above results to the case of Riemannian manifolds.
 
 For Riemannian manifolds, the $L^{p}$ boundedness of the Riesz transform of $-\Delta+V$ was discussed by many authors. We mention Meyer \cite{meyerP}, Bakry \cite{bakry} and Yosida \cite{yosidaN}. The most general answer was given by Sikora \cite{sikora1}. Let $M$ satisfying the doubling property $(D)$ and assume that the heat kernel verifies $\|p_t(x,.)\|_{2}\leq \frac{C}{\mu(B(x,\sqrt{t}))}$ for all $x\in M$ and $t>0$. Under these hypotheses, Sikora proved that if $V\in
L^{1}_{loc}(M)$, $V\geq 0$, then the Riesz transforms of $-\Delta+V$ are $L^p$ bounded for
$1<p\leq 2$ and of weak type $(1,1)$. 

Li \cite{li3} obtained boundedness results on Nilpotent Lie groups under the restriction $V\in RH_q$ and $q\geq \frac{D}{2}$,  $D$ being the dimension at infinity of $G$ (see \cite{coulhon7}).

Following the method of \cite{auscher3}, we obtain new results for $p>2$ on complete 
Riemannian manifolds satisfying the doubling property $(D)$, a Poincar\'e inequality $(P_{2})$ and taking $V$ in some $RH_q$. For manifolds of polynomial type we obtain additional results. This includes Nilpotent Lie groups.


Let us summarize the content of this paper. Let $M$ be a complete Riemannian manifold
 satisfying the doubling property $(D)$ and admitting a Poincar\'e inequality $(P_2)$. First we obtain the range of $p$ for the following maximal inequality valid for $u\in C^{\infty}_{0}(M)$:
 \begin{equation}\label{eq:Mp}
 \|\Delta u\|_{p}+ \|Vu\|_{p} \lesssim \|(-\Delta+V)u\|_{p}.
 \end{equation}
  The starting step is the following $L^1$ inequality for
$u\in C^{\infty}_{0}(M)$,
 \begin{equation}\label{max}
 \|\Delta u\|_{1}+ \|Vu\|_{1} \leq 3 \|(-\Delta+V)u\|_{1}
 \end{equation}
  which holds for any non-negative potential $V\in L^1_{loc}(M)$. This allows us to define $-\Delta+V$ as an operator on
$L^1(M)$ with domain $\mathcal{D}_{1}(\Delta) \cap \mathcal{D}_{1}(V)$. 

 For larger range of $p$, we assume that $V\in L^{p}_{loc}(M)$ and $-\Delta +V$ is a priori
defined on $C^{\infty}_{0}$. The validity of \eqref{eq:Mp} can be
obtained if one imposes for the
potential $V$ to be more regular:
\begin{thm}\label{th:1} Let $M$ be a complete Riemannian manifold satisfying $(D)$ and $(P_{2})$. Consider $V\in RH_{q}$ for some $1<q\le \infty$. Then there is 
$\epsilon>0$ depending only on $V$ such that \eqref{eq:Mp} holds for $1<p<q+\epsilon$.
 \end{thm}  
This new result for Riemannian manifolds is an extension of the one of Li \cite{li3} on Nilpotent Lie
groups settings obtained under the restriction $q\geq \frac{D}{2}$. 

The second purpose of our work is to establish some $L^p$ estimates for the
square root of $-\Delta+V$. 
Notice that we always have the identity
\begin{equation}\label{type2}
\|\,|\nabla  u|\,\|_2^2 + \|V^{\frac{1}{2}} u\|_2^2 = \|(- \Delta  + V)^{\frac{1}{2}}u \|_2^2, \quad
u\in C^{\infty}_{0}(M).
\end{equation}
The weak type $(1,1)$ inequality proved by Sikora \cite{sikora1} is satisfied under our hypotheses:
\begin{equation}
\label{eq:RT1}
\|\,|\nabla u|\, \|_{1,\infty} + \| V^{\frac{1}{2}}u \|_{1, \infty} \lesssim \| (- \Delta  +
V)^{\frac{1}{2}}u \|_{1}.
\end{equation}
 Interpolating (\ref{type2}) and (\ref{eq:RT1}), we obtain 
 \begin{equation}
\label{eq:RTp}
\|\,|\nabla u|\, \|_{p} + \| V^{\frac{1}{2}}u \|_{p} \lesssim \| (- \Delta  + V)^{\frac{1}{2}}u
\|_{p}
\end{equation} 
when  $1<p<2$ and  $u\in C_{0}^\infty(M)$. Here, $\|\ \|_{p,\infty}$ is  the norm in the Lorentz space
$L^{p,\infty}$ and $\lesssim$ is  the comparison in the sense of norms. 

It remains to find good assumptions on $V$ and $M$ to obtain (\ref{eq:RTp}) for some/all $2<p<\infty$. Recall before the following result
\begin{prop}\label{RTN} (\cite{auscher1}) Let $M$ be a complete Riemannian manifold satysfying $(D)$ and $(P_2)$. Then there exists $p_{0}>2$ such that the Riesz transform $T=\nabla (-\Delta)^{-\frac{1}{2}}$ is $L^{p}$ bounded for $1<p< p_{0}$.
\end{prop}
We now let $p_0=\sup\left\lbrace p\in ]2,\infty[; \nabla (-\Delta)^{-\frac{1}{2}} \textrm { is } L^p \; \textrm {bounded} \, \right\rbrace$. We obtain
the following theorem.

\begin{thm} \label{RH} Let $M$ be a complete Riemannian manifold. Let $V\in RH_q$ for some $q>1$ and $\epsilon>0$ such that $V\in RH_{q+\epsilon}$.
\begin{itemize} 
\item[1.] Assume that $M$ satisfies $(D)$ and $(P_2)$. Then for all $u\in C^{\infty}_0(M)$,
\begin{equation}\label{grad}
\|\,|\nabla
u|\,\|_{p}\lesssim\|(-\Delta+V)^{\frac{1}{2}}u\|_{p}\quad\textrm{
for } 1<p<\inf(p_0,2(q+\epsilon));
\end{equation}
\begin{equation}\label{V12}
\|V^{\frac{1}{2}}u\|_{p}\lesssim\|(-\Delta+V)^{\frac{1}{2}}u\|_{p}\quad\textrm{
for } 1<p<2(q+\epsilon).
\end{equation}
\item[2.] Assume that $M$ is of polynomial type and admits $(P_2)$. 
  Suppose that $D<p_{0}$, where $D$ is the dimension at infinity and that $\frac{D}{2}\leq q<\frac{p_{0}}{2}$. 
 \begin{itemize}
\item[a.] If $q<D$, then (\ref{grad}) holds for $1<p<\inf(q^{*}_{D}+\epsilon,p_{0})$,
 ($q^{*}_{D}=\frac{D q}{D-q})$.
\item[b.] If $q\geq D$, then (\ref{grad}) holds for $1<p<p_{0}$.
\end{itemize}
\end{itemize}
\end{thm}
Some remarks concerning this theorem:
\begin{itemize}
\item[1.]Note that point 1. is true without any additional assumption on the volume growth of balls other than $(D)$. Our assumption that $M$ is of polynomial type in point 2. --which is stronger than the doubling property (see section 2)-- is used only to improve the $L^p$ boundedness of $\nabla(-\Delta+V)^{-\frac{1}{2}}$ when $\frac{D}{2}<q<\frac{p_0}{2}$. We do not need it to prove $L^{p}$ estimates for $V^{\frac{1}{2}}(-\Delta+V)^{-\frac{1}{2}}$.
 \item[2.]If $q>\frac{p_0}{2}$ then we can replace $q$ in point 2. by any $q'<\frac{p_0}{2}$ since $V\in RH_{q'}$ (see Proposition \ref{CRO} in section 2).
 \item[3.]If $p_0\leq D$ and $q\geq \frac{D}{2}$, then (\ref{grad}) holds for $1<p<p_0$ and that is why we assumed $D<p_0$ in point 2..
 \item[4.] Finally the parameter $\epsilon$ depends on the self-improvement of the reverse H\"{o}lder condition (see Theorem \ref{CRO} in section 2).
 \end{itemize}

We establish also a converse theorem which is a crucial step in proving Theorem \ref{RH}.
\begin{thm}\label{RRT} Let $M$ be a  complete Riemannian manifold satisfying $(D)$ and $(P_l)$ for some $1\leq l< 2$. Consider $V\in RH_q$ for some $q>1$. Then
\begin{equation}\label{L1f}
\|(-\Delta+V)^{\frac{1}{2}}u\|_{l,\infty}\lesssim \|\,|\nabla
u|\,\|_{l}+\|V^{\frac{1}{2}}u\|_{l} \quad \textrm{ for every  }\;u\in C^{\infty}_{0}(M)
\end{equation}
and
\begin{equation}\label{Lp}
\|(-\Delta+V)^{\frac{1}{2}}u\|_{p}\lesssim \|\,|\nabla
u|\,\|_{p}+\|V^{\frac{1}{2}}u\|_{p} \quad \textrm{ for every  }\;u \in C^{\infty}_{0}(M)
\end{equation}
and  $l<p<2$.
\end{thm}

Using the interpolation result of \cite{badr2}, we remark that (\ref{Lp}) follows directly from (\ref{L1f}) and the  the $L^{2}$ estimate (\ref{type2}).
\begin{rem} The estimate (\ref{Lp}) always holds in the range $p>2$. This follows from the fact that (\ref{eq:RTp}) holds for $1<p\leq 2$ and that (\ref{eq:RTp}) for $p$ implies (\ref{Lp}) for $p'$, where $p'$ is the conjugate exponent of $p$.
\end{rem}
In the following corollaries we give examples of manifolds satisfying our hypotheses and to which we can apply the theorems above.
\begin{cor} Let $M$ be a  complete $n$-Riemannian manifold with non-negative Ricci
curvature. Then Theorem \ref{th:1},  part 1. of Theorem \ref{RH} and Theorem \ref{RRT} hold with $p_{0}=\infty$. Moreover, if $M$ satisfies the maximal volume growth $\mu(B)\geq cr^{n}$ for all balls $B$ of radius $r>0$ then part 2. of Theorem \ref{RH} also holds.
\end{cor}
\begin{proof} It suffices to note that in this case $M$ satisfies $(D)$ with $log_{2}C_{d}=n$, $(P_{1})$ --see Proposition \ref{DPR} below--, that the Riesz transform is $L^{p}$ bounded for $1<p<\infty$ \cite{bakry} and that $M$ has at most an Euclidean volume growth, that is $\mu(B)\leq Cr^{n}$ for any ball $B$ of radius $r>0$ --Theorem 3.9 in \cite{chavel}.
\end{proof}
\begin{cor} Let $C(N)=\mathbb{R}^{+}\times N$ be a conical manifold with compact basis $N$ of dimension $n-1\geq1$. Then Theorem \ref{th:1}, Theorem \ref{RRT} and Theorem \ref{RH} hold with $d=D=n$, $p_{0}=p_{0}(\lambda_1)>n$ where $\lambda_{1}$ is the first positive eigenvalue of the Laplacian on $N$.
\end{cor} 
\begin{proof}
Note that such a manifold is of polynomial type $n$ 
$$C^{-1}r^{n}\leq \mu(B)\leq Cr^{n}
$$
for all ball $B$ of $C(N)$ of radius $r>0$ (Proposition 1.3, \cite{li4}). $C(N)$ admits $(P_2)$ \cite{coulhon3}, and even $(P_1)$ using the methods in \cite{grigoryan2}. For the $L^p$ boundedness of the Riesz transform it was proved by Li \cite{li1} that $p_0=\infty$ when $\lambda_1\geq n-1$ and $p_0=\frac{n}{\frac{n}{2}-\sqrt{\lambda_{1}+(\frac{n-1}{2})^{2}}}>n$ when $\lambda_1<n-1$.
\end{proof}

Our main tools to prove these theorems are:
\begin{itemize}
\item the fact that $V$ belongs to a Reverse H\"{o}lder class;
\item an improved Fefferman-Phong inequality;
\item a Calder\'on-Zygmund decomposition;
\item reverse H\"older inequalities involving the weak solution of $-\Delta u+Vu=0$;
\item  complex interpolation;
\item the boundedness of the Riesz potential when $M$ satisfies  $\mu(B(x,r))\geq Cr^{\lambda}$ for all $r>0$.
\end{itemize}
Many arguments follow those of \cite{auscher3} --with additional technical problems due to the geometry of the Riemannian manifold-- but those for the Fefferman-Phong inequality require some sophistication. This Fefferman-Phong inequality with respect to balls is new even in the Euclidean case. In \cite{auscher3}, this inequality was proved with respect to cubes instead of balls which greatly simplifies the proof.

We end this introduction with a plan of the paper. In section 2, we recall the definitions of the doubling property, Poincar\'{e} inequality, reverse H\"{o}lder classes and homogeneous Sobolev spaces associated to a potential $V$. Section 3 is devoted to define the Schr\"{o}dinger operator. In section 4 we give the principal tools to prove the theorems mentioned above. We establish an improved Fefferman-Phong inequality, make a Calder\'{o}n-Zygmund decomposition, give estimates for positive subharmonic functions. We prove Theorem \ref{th:1} in section 5. We handle the proof of Theorem \ref{RH}, part 1. in section 6. Section 7 is concerned with the proof of Theorem \ref{RRT}. In section 8, we give different estimates for the weak solution of $-\Delta u +V u=0$ and complete the proof of item 2. of Theorem \ref{RH}. Finally, in section 9, we apply our result on Lie groups with polynomial growth. 
 \\

  \textit{Acknowledgements.} The two authors would like to thank their Ph.D advisor P. Auscher for proposing this joint work and for the useful discussions and advice on the topic of the paper.
 \section{Preliminaries}
  Let $M$ denote a complete non-compact Riemannian manifold. We write $\rho$ for the geodesic distance, $\mu$ for the Riemannian measure on $M$, $\nabla$ for the Riemannian gradient, $\Delta$ for the Laplace-Beltrami operator, $|\cdot|$ for the length on the tangent space (forgetting the subscript $x$ for simplicity) and $\|\cdot\|_{p}$ for the norm on $ L^{p}(M,\mu)$, $1 \leq p\leq +\infty.$
\subsection{The doubling property and Poincar\'{e} inequality}
\begin{dfn}[Doubling property] Let $M$ be a Riemannian manifold. Denote by $B(x, r)$ the open ball of center $x\in M $ and radius $r>0$. One says that $M$ satisfies the doubling property $(D)$ if there exists a constant $C_{d}>0$, such that for all $x\in M,\, r>0 $ we have
\begin{equation*}\tag{$D$}
\mu(B(x,2r))\leq C_{d} \mu(B(x,r)).
\end{equation*}
\end{dfn}
\begin{lem} Let $M$ be a Riemannian manifold satisfying $(D)$ and let $s=log_{2}C_{d}$. Then for all $x,\,y\in M$ and $\theta\geq 1$
\begin{equation}\label{teta}
\mu(B(x,\theta R))\leq C\theta^{s}\mu(B(x,R))
\end{equation}
and
\begin{equation}\label{y}
\mu(B(y, R))\leq C(1+\frac{d(x,y)}{R})^{s}\mu(B(x,R)).
\end{equation}
\end{lem} 
We have also the following lemma:
 \begin{lem}\label{CD} Let $M$ be a Riemannian manifold satisfying $(D)$. Then for $x_{0}\in M$, $r_{0}>0$, we have
 $$
 \frac{\mu(B(x,r))}{\mu(B(x_{0},r_{0}))}\geq 4^{-s} (\frac{r}{r_{0}})^{s}
 $$
 whenever $x\in B(x_{0},r_{0})$ and $r\leq r_{0}$.
 \end{lem}
 \begin{thm}[Maximal theorem]\label{MIT} (\cite{coifman2})
Let $M$ be a Riemannian manifold satisfying $(D)$. Denote by $\mathcal{M}$ the uncentered Hardy-Littlewood maximal function over open balls of $X$ defined by
 $$
 \mathcal{M}f(x)=\underset{B:x\in B}{\sup}|f|_{B}
 $$ 
 where $ \displaystyle f_{E}:=\aver{E}f d\mu:=\frac{1}{\mu(E)}\int_{E}f d\mu.$
Then

\begin{itemize}
\item[1.] $\mu(\left\lbrace x:\,\mathcal{M}f(x)>\lambda\right\rbrace)\leq \frac{C}{\lambda}\int_{M}|f| d\mu$ for every $\lambda>0$;
\item[2.] $\|\mathcal{M}f\|_{p}\leq C_{p} \|f\|_{p}$, for $1<p\leq\infty$.
\end{itemize}
\end{thm}
\begin{dfn}A Riemannian manifold $M$ is of polynomial type if there is $c,\,C>0$ such that
 \begin{equation} \tag{$LU_l$}
c^{-1}r^{d} \leq\mu(B(x,r))\leq cr^{d}
 \end{equation}
  for all $x\in M$ and $r\leq 1$ and
   \begin{equation} \tag{$LU_\infty$}
C^{-1}r^{D} \leq\mu(B(x,r))\leq Cr^{D}
 \end{equation}
  for all $x\in M$ and $r\geq 1$.
  \end{dfn}
  We call $d$ the local dimension and $D$ the dimension at infinity.
  Note that if $M$ is of polynomial type then it satisfies $(D)$ with $s=\max(d,D)$. Moreover, for every $\lambda\in[\min(d,D), \max(d,D)]$,
   \begin{equation} \tag{$L_\lambda$}
\mu(B(x,r))\geq cr^{\lambda}
 \end{equation}
 for all $x\in M$ and $r>0$.
 \\
 \begin{dfn}[Poincar\'{e} inequality] Let $M$ be a complete Riemannian manifold, $1\leq l<\infty$. We say that $M$ admits a Poincar\'{e} inequality $(P_{l})$ if there exists  a constant $C>0$ such that, for every function $f\in C^{\infty}_{0}(M)$, and every ball $B$ of $M$ of radius $r>0$, we have
\begin{equation*}\tag{$P_{l}$}
\left(\aver{B}|f-f_{B}|^{l} d\mu\right)^{\frac{1}{l}} \leq C r\left(\aver{B}|\nabla f|^{l}d\mu\right)^{\frac{1}{l}}.
\end{equation*}
\end{dfn}
\begin{rem} Note that if $(P_{l})$ holds for all $f\in C_{0}^{\infty}$, then it holds for all $f\in W_{p,loc}^1$ for $p\geq l$ (see \cite{hajlasz4}, \cite{keith2}).
\end{rem}
The following result from Keith-Zhong \cite{keith2} improves the exponent in the Poincar\'e inequality.
\begin{lem}\label{kz} Let $(X,d,\mu)$ be a complete metric-measure space satisfying $(D)$
and admitting a Poincar\'{e} inequality $(P_{l})$, for  some $1< l<\infty$.
Then there exists $\epsilon >0$ such that $(X,d,\mu)$ admits
$(P_{p})$ for every $p>l-\epsilon$. 
\end{lem}
\begin{prop}\label{DPR} Let $M$ be a complete Riemannian manifold $M$  with non-negative Ricci curvature. Then $M$ satisfies $(D)$ (with $C_{d}=2^n$) and admits a Poincar\'{e} inequality $(P_{1})$.
  \end{prop}
  \begin{proof}
  Indeed if the Ricci curvature is non-negative that is there exists
  $a>0$ such that 
  $R_{ic} \geq -a^{2}g$, 
  a  result by Gromov \cite{cheeger2}
  shows that
  $$\mu(B(x,2r))\leq 2^{n}\mu (B(x,r)) \,\textrm{ for all}\, x\in M, \; r>0.$$
  Here $n$ means the topologic dimension.\\
  On the other hand, Buser's inequality \cite{buser} gives us
  $$\int_{B}|u-u_{B}|\, d\mu \leq c(n)  r \int_{B}|\nabla u| \,
  d\mu.
  $$
 Thus we get $(D)$ and $(P_{1})$ (see also \cite{saloff1}).
 \end{proof}
 \subsection{Reverse H\"{o}lder classes}
\begin{dfn} Let $M$ be a Riemannian manifold. A weight $w$ is a non-negative locally integrable function on $M$. The reverse H\"{o}lder classes are defined in the following way: $w\in RH_{q},\,1<q<\infty$, if
\begin{itemize}
\item[1.] $wd\mu$ is a doubling measure;
\item[2.] there exists a constant $C$ such that for every ball $B\subset M$
\begin{equation}\label{rhq}
\left(\aver{B}w^{q}d\mu\right)^{\frac{1}{q}}\leq C\aver{B}wd\mu.
\end{equation}
The endpoint $q=\infty$ is given by the condition: $w\in RH_{\infty}$ whenever, $wd\mu$ is doubling and for any ball $B$,
\begin{equation}\label {rhi}
w(x)\leq C\aver{B}w \quad \textrm{for } \mu-a.e. \;x\in B.
\end{equation}
\end{itemize}
\end{dfn}
On $\mathbb{R}^n$, the condition $wd\mu$ doubling is superfluous. It could be the same on a Riemannian manifold.
\begin{prop}\label{CRO}(\cite{stromberg}, \cite{garcia})
\begin{itemize}
\item[1.] $RH_{\infty}\subset RH_{q}\subset RH_{p}$ for $1<p\leq q\leq \infty$.
\item[2.] If $w\in RH_{q}$, $1<q<\infty$, then there exists $q<p<\infty$ such that $w\in RH_{p}$.
\item[3.] We say that $w\in A_p$ for $1<p<\infty$ if there is a constant $C$ such that for every ball $B\subset M$
$$\left(\aver{B}wd\mu\right)\left(\aver{B}w^{\frac{1}{1-p}}d\mu\right)^{p-1}\leq C.$$ 
For $p=1$, $w\in A_1$ if there is a constant $C$ such that for every ball $B\subset M$
$$
\aver{B}wd\mu\leq Cw(y)\quad \textrm {for }\;\mu-a.e. y\in B.
$$  
 We let $A_{\infty}=\bigcup_{1\leq p<\infty}A_p$. Then $A_{\infty}=\bigcup_{1<q\leq \infty} RH_{q}$. 
\end{itemize}
\end{prop}
\begin{prop}\label{CR}(see section 11 in \cite{auscher3}, \cite{johnson})
Let $V$ be a non-negative measurable function. Then the following properties are equivalent:
\begin{itemize}
\item[1.] $V\in A_{\infty}$.
\item[2.] For all $r\in]0,1[,\,V^{r}\in RH_{\frac{1}{r}}$.
\item[3.] There exists $r\in ]0,1[,\,V^{r}\in RH_{\frac{1}{r}}$.
\end{itemize}
\end{prop}
We end this subsection with the following lemma:
\begin{lem}\label{buck2} Let $G$ be an open subset of an homogeneous space $(X,d,\mu)$ and let $\mathcal{F}(G)$ be the set of metric balls contained in $G$. Suppose that for some $0<q<p$ and non-negative $f\in L^{p}_{loc}$, there is a constant $A>1$ and $1\leq \sigma_{0}\leq \sigma'_{0}$ such that
$$
\left(\aver{B}f^{p}d\mu\right)^{\frac{1}{p}}\leq A\left(\aver{\sigma_{0}B}f^{q}d\mu\right)^{\frac{1}{q}}\quad \forall B:\,\sigma'_{0}B\in \mathcal{F}(G).
$$
Then for any $0<r<q$ and $1<\sigma\leq \sigma'<\sigma'_{0}$, there exists a constant $A'>1$ such that
$$
\left(\aver{B}f^{p}d\mu\right)^{\frac{1}{p}}\leq A'\left(\aver{\sigma B}f^{r}d\mu\right)^{\frac{1}{r}} \quad \forall B:\, \sigma'B\in \mathcal{F}(G).
$$
\end{lem}

\subsection{Homogeneous Sobolev spaces associated to a weight $V$}
\begin{dfn}(\cite{badr2})
Let $M$ be a Riemannian manifold, $V\in A_{\infty}$. Consider  for $1\leq p<\infty$, the vector space $\dot{W}_{p,V}^{1}$ of 
 distributions $ f $ such that $|\nabla f|$ and $ Vf \in L^{p}$. 
 It is well known that the elements of $\dot{W}_{p,V}^{1}$ are in  $L^{p}_{loc} $. We equip $\dot{W}_{p,V}^{1}$  with the semi norm 
$$
\|f\|_{\dot{W}_{p,V}^{1}}=\|\,|\nabla f|\,\|_{p}+ \|Vf\|_{p}.
$$
\end{dfn}
In fact, this expression is a norm since $V\in A_{\infty}$ yields $V>0 \; \mu-a.e$.
\begin{dfn} We denote $\dot{W}_{\infty,V}^{1}$ the space of all Lipschitz functions $f$ on $M$ with $\|Vf\|_{\infty}<\infty$.
\end{dfn}
\begin{prop}(\cite{badr2}) \label{BS}Assume that $M$ satisfies $(D)$ and admits a Poincar\'{e} inequality $(P_{s})$ for some $1\leq s<\infty$ and that $V\in A_{\infty}$. Then, for $s\leq p\leq \infty$, $\dot{W}_{p,V}^{1}$ is a Banach space.
\end{prop} 
\begin{prop} \label{RW}Under the same hypotheses as in Proposition \ref{BS}, the Sobolev space $\dot{W}_{p,V}^{1}$ is reflexive for $s\leq p<\infty$.
\end{prop}
\begin{proof} The Banach space $\dot{W}_{p,V}^{1}$ is isometric to a closed subspace of $L^{p}(M,\mathbb{R}\times T^{*}M)$ which is reflexive. The isometry is given by the linear operator $T: \dot{W}_{p,V}^{1}\rightarrow L^{p}(M,\mathbb{R}\times T^{*}M)$ such that $Tf=(Vf,\nabla f)$  by definition of the norm of $\dot{W}_{p,V}^{1}$ and Proposition \ref{BS}.
\end{proof}
 \begin{thm} \label{rt8} (\cite{badr2}) Let $M$ be a complete Riemannian manifold satisfying $(D)$. Let $V\in RH_{q}$  for some $1< q\leq\infty$ and assume that $M$ admits a Poincar\'{e} inequality $(P_{l})$ for some $1\leq l<q$. Then, for $1\leq p_{1}<p<p_{2}\leq q$, with $p>l$, $\dot{W}_{p,V}^{1}$ is a real interpolation space between $\dot{W}_{p_1,V}^{1}$ and $\dot{W}_{p_2,V}^{1}$.
 \end{thm}

\section{Definition of Schr\"odinger operator}
Let $V$ be a non-negative, locally integrable function on M. Consider the
sesquilinear form $$\mathcal{Q}(u,v) = \int_{M} (\nabla u \cdot \overline{\nabla v }
+ V u \, \overline v)d\mu $$
 with domain $$\mathcal{V}=\mathcal{D}(\mathcal{Q}) =W_{2,V^{\frac{1}{2}}}^{1}= \{f\in L^2(M)\, ; \,
|\nabla f| \ \& \ V^{\frac{1}{2}}f \in L^2(M)\}$$
equipped with the norm
$$
\|f\|_{\mathcal{V}}=(\|f\|_{2}^{2}+\|\nabla f\|_{2}^{2}+\|V^{\frac{1}{2}}f\|_{2}^{2})^{\frac{1}{2}}.
$$
 Clearly $\mathcal{Q}(.,.) $ is a positive, symmetric closed
form. It follows that there exists a unique positive self-adjoint operator, which we
 call
$H=-\Delta +V$ , such that
$$
\langle Hu, v\rangle = \mathcal{Q}(u,v) \quad \forall\, u \in \mathcal{D}(H), 
\;\forall\, v\in \mathcal{V}.
$$
When $V=0$, $H=-\Delta$ is the Laplace-Beltrami operator.
Note that $C^{\infty}_{0} (M)$ is dense in $\mathcal{V}$ (see the Appendix in \cite{badr2}).

The Beurling-Deny theory holds on $M$, which means that
$\epsilon(H+\epsilon)^{^-1}$ is a positivity-preserving contraction on  $L^p(M)$  
for all $1\leq p \leq \infty$ and $\epsilon>0$. 
Moreover, if $V'\in L^{1}_{loc}(M)$ such that $0\le
V'\le V$ and $H'$ is the corresponding operator then one has 
for any $\epsilon>0$ and  for any $f \in L^p$, $1\le p \le \infty$, $f\ge 0$
$$0 \le(H+\epsilon)^{-1} f \leq (H' + \epsilon)^{-1} f .$$  
It is equivalent to a pointwise comparision of the kernels of resolvents. In
particular, if  $V$ is bounded  from below by some positive constant $\epsilon>0$,
then $H^{-1}$ is bounded on $L^p$ for $1\le p \le \infty$  and is
dominated by  $ (-\Delta + \epsilon)^{-1}$ (see Ouhabaz \cite{ouhabaz}).

Let $\dot{\mathcal{V}}$ be the closure of 
$C^{\infty}_{0}(M)$ under the semi-norm 
$$
\|f\|_{\dot{\mathcal{V}}}= \big( \|\,|\nabla f|\,\|_{2}^2 + \|V^{\frac{1}{2}}
f\|_{2}^2\big)^{\frac{1}{2}}.
$$
Assume that $M$ satisfies $(D)$ and $(P_{2})$.
By  Fefferman-Phong inequality --Lemma \ref{lem:FP} in section 4 below--, there is a continuous inclusion $\dot{\mathcal{V}} \subset L^{2}_{loc}$ if $V$ is not identically 0, which is assumed from now on, 
hence, this is a norm.    
Let $f\in \dot{\mathcal{V}}'$. Then, there exists a unique $u\in
\dot{\mathcal{V}}$ such that 
\begin{equation}
\label{eq:Htilde}
\int_{M} \nabla u \cdot {\nabla \overline v } + V u\,  \overline v = \langle f,
v \rangle \quad \forall\, v\in C^{\infty}_{0}(M).
\end{equation}
In particular, $-\Delta u + Vu = f$ holds in the distributional sense.
We can obtain $u$ for a nice $f$ by the next lemma. 

\begin{lem}\label{lem:approx}Assume that $M$ satisfies $(D)$ and $(P_{2})$.
Consider $f\in C^{\infty}_{0}(M) \cap L^2(M)$. For $\epsilon>0$, let
$u_{\epsilon}= (H+\epsilon)^{-1}f \in \mathcal{D}(H)$. Then 
$(u_{\epsilon})$ is a bounded sequence in $\dot{\mathcal{V}}$ which
converges strongly to $H^{-1} f$.
\end{lem}
\begin{proof}
The proof is analogous to the  proof of Lemma 3.1 in \cite{auscher3}.
\end{proof}
\begin{rem}Assume that $M$ satisfies $(D)$ and $(P_{2})$.
The continuity of the inclusion $\dot {\mathcal{V}} \subset L^2_{loc}(M)$ has two further consequences. First, we have that $L^2_{comp}(M)$, the space of compactly supported $L^2$ functions on $M$, is continuously contained in 
$\dot{\mathcal{V'}} \cap L^{2}(M)$.  Second, $(u_{\epsilon})$ has a subsequence converging to $u$ almost everywhere.
\end{rem}

 Finally as $H$ is self-adjoint, it has a unique square root which we denote $H^{\frac{1}{2}}$. $H^{\frac{1}{2}}$ is
defined as the unique maximal-accretive operator such that $H^{\frac{1}{2}}H^{\frac{1}{2}} =H$.
 We have that $H^{\frac{1}{2}}$ is self-adjoint with domain $\mathcal{V}$ and for all $u\in
C_{0}^\infty(M)$, $\|H^{\frac{1}{2}}u\|_{2}^2= \|\,|\nabla u|\,\|_{2}^2+\|V^{\frac{1}{2}}
u\|_{2}^2$.  This allows us to extend 
$H^{\frac{1}{2}}$ from  $\dot{\mathcal{V}}$ into $L^2(M)$. If $S$ denotes this extension,
then we have $S^{\star}S=H$ where $S^{\star}\colon L^2(M) \to \dot{\mathcal{V}}'$
is the adjoint of $S$.

\section{Principal tools}
We gather in these section the main tools that we need to prove our theorems. Some of them are of independent interest.
\subsection{An improved Fefferman-Phong inequality}
\begin{lem}\label{lem:FP} Let $M$ be a complete Riemannian manifold satisfying $(D)$. Let $w\in A_{\infty}$ and $1\leq p <\infty$. Assume that $M$ admits also a Poincar\'{e} inequality $(P_{p})$. Then there is a constant $C>0$ depending only on the $A_{\infty}$ constant of $w$, $p$ and the constants in $(D),\,(P_{p})$, such that for every ball $B$ of radius $R>0$ and $u \in W_{p,loc}^{1}$
\begin{equation}\label{eq:FP}
\int_{B} (|\nabla u|^p + w|u|^p )d\mu\geq \frac{C m_{\beta}({R^{p} w_{B}})}{ R^{p}} \int_B|u|^pd\mu 
\end{equation}
where $m_{\beta}(x)  = x$ for $x\leq 1$ and $m_{\beta}(x) = x^\beta$ for $x\geq 1$.
\end{lem}
\begin{proof}
 Since $M$ admits a $(P_{p})$ Poincar\'{e} inequality, we have
$$
\int_{B}|\nabla u|^{p}d\mu \geq \frac{C}{R^{p}\mu(B)} \int_{B}\int_{B} |u(x)-u(y)|^{p}d\mu(x)d\mu(y).
$$
 This and 
$$
\int_{B}w|u|^{p}d\mu =\frac{1}{\mu(B)}\int_{B}\int_{B} w(x) |u(x)|^{p} d\mu(x)d\mu(y)
$$
lead easily to
$$
\int_{B}(|\nabla u|^{p}+w|u|^{p})d\mu \geq [\min(CR^{-p},w)]_{B}\int_{B}|u|^{p}d\mu.
$$

Now we use that $w\in A_{\infty}$. There exists $\varepsilon>0$, independent of $B$, such that $\,E=\left\lbrace x\in B: w(x)>\varepsilon w_{B}\right\rbrace $ satisfies $\mu(E)>\frac{1}{2}\mu(B)$. Hence 
$$
[\min(CR^{-p},w)]_{B}\geq \frac{1}{2} \min(CR^{-p},\varepsilon w_{B})\geq C \min(R^{-p},w_{B}).
$$
This proves the desired inequality 
when $R^{p}w_{B}\leq1$.

Assume now $R^{p}w_{B}>1$. We say that a ball $B$ of radius $R$ is of type 1 if $R^{p}w_{B}<1$ and of type 2 if not. Take $\delta,\,\epsilon>0$ such that $2\delta<\epsilon<1$. We consider a maximal covering of $(1-\epsilon)B$ by balls $(B_{i}^{1})_{i}:=(B(x_{i}^{1},\delta R))_{i}$ such that the balls $\frac{1}{2}B_{i}^{1}$ are pairwise disjoint. By $(D)$ there exists $N$ independent of $\delta$ and $R$ such that $\sum_{i\in I}\ind_{B_{i}^{1}}\leq N$. Since $2\delta<\epsilon$, we have $B_{i}^{1}\subset B$ for all $i\in I$.
Denote $G_{1}$ the union of all balls $B_{i}^{1}$ of type 1 and $\widetilde{G}_{1}=\{ x\in M:\,d(x,G_{1})\leq \epsilon \delta R\}$. Set $\widetilde{E}_{1}=(1-\epsilon \delta)B-\widetilde{G}_{1}$. This time we consider a maximal covering of $\widetilde{E}_{1}$ by balls $(B_{i}^{2})_{i}:=(B(x_{i}^{2},\delta^{2} R))_{i}$ such that the balls $\frac{1}{2}B_{i}^{2}$ are pairwise disjoint. Therefore with the same $N$ one has $\sum_{i\in I}\ind_{B_{i}^{2}}\leq N$. Let $G_{2}$ be the union of all balls $B_{i}^{2}$ of type 1 and $\widetilde{G}_{2}=\{ x\in M:\,d(x,G_{1}\cup G_{2})\leq \epsilon \delta^{2} R \}$, $\widetilde{E}_{2}=(1-\epsilon \delta^{2})B-\widetilde{G}_{1}$. We iterate this process. Note that the $G_{j}$'s are pairwise disjoint (from $2\delta<\epsilon$). We claim  then that $\mu(B-\bigcup_{j}G_{j})=0$. Indeed, for almost $x\in B$, $w_{B'}$ converges to $w(x)$ whenever $r(B')\rightarrow 0$ and $x\in B'$. Take such an $x$ and assume that $x\notin \bigcup_{j}G_{j}$. Then, for every $j$ there exists $x_{k}^{j}$ such that $x\in B(x_{k}^{j},\delta^{j}R)$ and $(\delta^{j}R)^{p}w_{B(x_{k}^{j},\delta^{j}R)}\geq1$. This is a contradiction since $(\delta^{j}R)^{p}w_{B(x_{k}^{j},\delta^{j}R)}\rightarrow 0$ when $j\rightarrow \infty$. Note also that there exists $0<A<1$ such that for all $j,\,k$ and ball $B_{k}^{j}$ of type 1, 
\begin{equation}\label{wD}
(\delta^{j}R)^{p}w_{B_{k}^{j}}> A.
\end{equation}
 Indeed, let $B_{k}^{j}$ be of type 1. There exists $B_{l}^{j-1}$ such that $x_{k}^{j}\in B_{l}^{j-1}$ and $B_l^{j-1}$ must be of type 2 because $x_{k}^{j}\notin G_{j-1}$. Hence $B_{k}^{j}\subset B(x_{l}^{j-1}, \delta^{j}(1+\delta^{-1})R)$. Since $wd\mu$ is doubling, we get
\begin{align*}
w(B_{l}^{j-1})&\leq w\left(B(x_{l}^{j-1},\delta^{j}(1+\delta^{-1})R)\right)
\\
&\leq C'(1+\delta^{-1})^{s'}w\left(B(x_{l}^{j-1},\delta^{j}R)\right)
\\
&\leq C'^{2}(1+\delta^{-1})^{s'}(1+\frac{d(x_{l}^{j-1},x_{k}^{j})}{\delta^{j}R})^{s'}w(B_{k}^{j})
\\
&\leq C'^{2}(1+\delta^{-1})^{2s'}w(B_{k}^{j})
\end{align*}
where $s'=log_{2}C'$ and $C'$ is the doubling constant of $wd\mu$.
 On the other hand, since $d\mu$ is doubling
 \begin{align*}
 \mu(B_{l}^{j-1})&\geq C^{-1}(1+\delta)^{-s}\mu(B(x_{l}^{j-1},\delta^{j-1}(1+\delta)R))
 \\
 &\geq C^{-1} (1+\delta)^{-s}\mu(B_{k}^{j}).
 \end{align*}
 Since $B_{l}^{j-1}$ is of type 2, we obtain 
 \begin{align*}
 (\delta^{j}R)^{p}w_{B_{k}^{j}}&\geq C'^{-2}C^{-1}(1+\delta^{-1})^{-2s'}(1+\delta)^{-s}\delta^{p}(\delta^{j-1}R)^{p}w(B_{l}^{j-1})
 \\
 &>C'^{-2}C^{-1}(1+\delta^{-1})^{-2s'}(1+\delta)^{-s}\delta^{p}.
 \end{align*}
  Thus we get (\ref{wD}) with $A=C'^{-2}C^{-1}(1+\delta^{-1})^{-2s'}(1+\delta)^{-s}\delta^{p}$.
 From all these facts we deduce that
\begin{align*}
\int_{B}(|\nabla u|^{p}+w|u|^{p})d\mu&\geq \frac{1}{N}\sum_{j,\,k: B_{k}^{j} \textrm{of type 1 }\,}\int_{B_{k}^{j}}(|\nabla u|^{p}+w|u|^{p})d\mu
\\
&\geq C\frac{1}{N}\sum_{j,\,k: B_{k}^{j} \textrm{of type 1 }\,}\min((\delta^{j}R)^{-p},w_{B_{k}^{j}})\int_{B_{k}^{j}}|u|^{p}d\mu
\\
&\geq \frac{C}{N}A\sum_{j,\,k: B_{k}^{j} \textrm{of type 1 }}\,(\delta^{j}R)^{-p}\int_{B_{k}^{j}}|u|^{p}d\mu
\\
&\geq \frac{C}{N} A\min_{j}\left( \frac {R} {\delta^{j} R} \right)^p R^{-p} \int_B|u|^pd\mu.
\end{align*}
We used Fefferman-Phong inequality in the second estimate, (\ref{wD}) in the penultimate one, and that the $B_k^j$ of type 1 cover $B$ up to a $\mu-$ null set in the last one.
It remains to estimate $ \min_{j} \left( \frac {R} {R_{j}} \right)^p$ from below with $R_j=\delta^jR$.
Let $1\leq \alpha <\infty$ be such that $w\in A_{\alpha}$ --the Muckenhoupt class--. Then for any ball $B$ and measurable subset $E$ of $B$ we have

$$
 \left(\frac{w_{E}}{w_{B}}\right)\geq C\left(\frac{\mu(E)}{\mu(B)}\right)^{\alpha-1}.
$$
Applying this to $E=B_{k}^{j}$ and $B$ we obtain
\begin{align*}
\left(\frac{R}{R_{j}}\right)^{p}&= \frac{R^{p}w_{B}}{R_{j}^{p}w_{B_{k}^{j}}} \frac{w_{B_{k}^{j}}}{w_{B}}
\\
&\geq R^{p} w_{B} \frac{w_{B_{k}^{j}}}{w_{B}}
\\
&\geq C R^{p}w_{B}\left(\frac{\mu(B_{k}^{j})}{\mu(B)}\right)^{\alpha-1}
\\
&\geq CR^{p}w_{B}\left(\frac{R_{j}}{R}\right)^{s(\alpha-1)}
\end{align*}
where we used Lemma \ref{CD}. 
This yields $\min_{j}(\frac{R}{R_{j}})^{p}\geq C(R^{p}w_{B})^{\beta}$ with $\beta=\frac{p}{p+s(\alpha-1)}$. The lemma is proved.
\end{proof}
\subsection{Calder\'{o}n-Zygmund decomposition}
We now proceed to establish the following Calder\'{o}n-Zygmund decomposition:
\begin{prop}\label{CZ} Let $M$ be a complete Riemannian manifold  satisfying $(D)$ and $(P_{l})$ for some $1\leq l<2$. Let $l\leq p<2$, $V\in A_{\infty}$, $f\in \dot{W}_{p,V^{\frac{1}{2}}}^{1}$ and $\alpha>0$. Then, one can find a collection of balls $(B_i)$, functions $g$ and $b_i$  satisfying the following properties:
\begin{equation}\label{eqcsds1}
f= g+\sum_i b_i  \end{equation}
\begin{equation}\label{eqcsds2}
\|\,|\nabla g|\,\|_2 + \|V^{\frac{1}{2}}g\|_{2} \leq C\alpha^{1-\frac{p}{2}}(\|\,|\nabla f|\,\|_{p} +\|V^{\frac{1}{2}}f\|_{p})^{\frac{1}{2}}, 
\end{equation}
\begin{equation}\label{eqcsds3}
\supp\,b_i  \subset{B_i}\ \text{and} \ \int_{B_i} (|\nabla b_i|^{l} + |V^{\frac{1}{2}}b_{i}|^{l}+ R_{i}^{-l}| b_i|^l)d\mu \leq C\alpha^l \mu(B_i), \end{equation}
 \begin{equation}\label{eqcsds4}
\sum_i \mu(B_i) \leq C\alpha^{-p} \int_{M} (|\nabla f|^p +  | V^{\frac{1}{2}}f |^p)d\mu, \end{equation}
\begin{equation}\label{eqcsds5}
\sum_i \ind_{B_i} \leq N, \end{equation}
where 
 $N$ depends only on the doubling constant, and $C$ on the doubling constant, $p,\,l$ and the $A_{\infty}$ constant of $V$. Here, $R_{i}$ denotes the radius of $B_{i}$ and gradients are taken in the distributional sense on $M$.
\end{prop}
\begin{rem}
The function $g$ is Lipschitz with Lipschitz constant controlled by $C\alpha$.
\end{rem} 
\begin{proof}[Proof of Proposition \ref{CZ}] Let $\Omega$ be the open set $ \{x \in M; \mathcal{M}(|\nabla f|^l+ | V^{\frac{1}{2}}f |^l)(x) >\alpha^l\}$.
If $\Omega$ is empty, then set $g=f$ and $b_{i}=0$. Otherwise, the maximal theorem --Theorem \ref{MIT}-- yields
\begin{equation*} \label{mO}
\mu(\Omega) \leq C\alpha^{-p} \int_{M} (|\nabla f|^p+| V^{\frac{1}{2}}f |^p)d\mu.  \end{equation*} 
  In particular $\Omega\neq M$ as $\mu(M)=\infty$. Let $F$ be the complement of $\Omega$. Since $\Omega$ is an open set distinct of $M$, let $(\underline{B_{i}})$ be a Whitney decomposition of $\Omega$ (\cite{coifman1}). That is, the balls $\underline{B_i}$ are pairwise disjoint and there is two constants $C_{2}>C_{1}>1$, depending only
on the metric, such that
\begin{itemize}
\item[1.] $\Omega=\bigcup_{i}B_{i}$ with $B_{i}=
C_{1}\underline{B_{i}}$ are contained in $\Omega$ and the balls $(B_{i})_{i}$ have the bounded overlap property;
\item[2.] $r_{i}=r(B_{i})=\frac{1}{2}d(x_{i},F)$ and $x_{i}$ is 
the center of $B_{i}$;
\item[3.] each ball $\overline{B_{i}}=C_{2}\underline{B_{i}}$ intersects $F$ ($C_{2}=4C_{1}$ works).
\end{itemize}
For $x\in \Omega$, denote $I_{x}=\left\lbrace i:x\in B_{i}\right\rbrace$. By the bounded overlap property of the balls $B_{i}$, we have that $\sharp I_{x} \leq N$. Fixing $j\in I_{x}$ and using the properties of the $B_{i}$'s, we easily see that $\frac{1}{3}r_{i}\leq r_{j}\leq 3r_{i}$ for all $i\in I_{x}$. In particular, $B_{i}\subset 7B_{j}$ for all $i\in I_{x}$.

Condition \eqref{eqcsds5} is nothing but the bounded overlap property of the $B_{i}$'s  and \eqref{eqcsds4} follows from \eqref{eqcsds5} and (\ref{mO}).
 We remark that since $V\in A_{\infty}$, Proposition \ref{CR} yields $V^{\frac{l}{2}} \in A_{\infty}$. Applying Lemma \ref{lem:FP}, we obtain
\begin{equation}\label{t8}
\int_{B_{i}} (|\nabla f|^l+ |V^{\frac{1}{2}}f|^l)d\mu \geq C \min ((V^{\frac{l}{2}})_{B_{i}} , R_{i}^{-l} )\int_{B_{i}}|f|^ld\mu. 
\end{equation}
We declare $B_{i}$ of type 1 if 
  $(V^{\frac{l}{2}})_{B_{i}} \geq R_{i}^{-l}$
 and of type 2 if $(V^{\frac{l}{2}})_{B_{i}} <  R_{i}^{-l}$.

Let us now define the functions $b_i$. Let $(\chi_i)$ be a partition of unity on $\Omega$
associated to the covering $(\underline{B_i})$ so that for each $i$, $\chi_i$ is a $C^1$ function supported in $B_i$ with $\|\chi_i\|_\infty +
R_i \|\,|\nabla \chi_i|\,\|_\infty \leq C$.   
Set
$$
b_i = \begin{cases} f\chi_i, & \mathrm{if}\ B_{i}\ \mathrm{is\ of\ type\ 1}, 
\\
 (f-f_{B_{i}})\chi_i, & \mathrm{if}\ B_{i}\ \mathrm{is \ of \ type\  2}.
\end{cases}
$$
 If $B_{i}$ is of type 2, then it is a direct consequence of the Poincar\'e inequality $(P_{l})$ that 
$$
\int_{B_i} (|\nabla b_i|^l +  R_{i}^{-l}| b_i|^l)d\mu \leq C \int_{B_{i}}  |\nabla f|^ld\mu.
$$
As $\int_{\overline{B_{i}}} |\nabla f|^ld\mu \leq \alpha^l  \mu(\overline{B_i})$ we get the desired inequality in  \eqref{eqcsds3}.
For $V^{\frac{1}{2}}b_{i}$ we have
\begin{align*}
\int_{B_{i}} |V^{\frac{1}{2}}b_{i}|^{l}d\mu &=\int_{B_{i}} |V^{\frac{1}{2}}(f-f_{B_{i}})\chi_{i}|^{l}d\mu
\\
&\leq C\left(\int_{B_{i}} |V^{\frac{1}{2}}f|^{l}d\mu+\int_{B_{i}}|V^{\frac{1}{2}}f_{B_{i}}|^{l}d\mu\right)
\\
&\leq C\left((|V^{\frac{1}{2}}f|^{l})_{B_{i}}\mu(B_{i})+C (V^{\frac{l}{2}})_{B_{i}}(|f|^{l})_{B_{i}}\mu(B_{i})\right)
\\
&\leq C \left(\alpha^{l} \mu(B_{i})+\left(|\nabla f|^{l}+|V^{\frac{1}{2}}f|^{l}\right)_{B_{i}}\mu(B_{i})\right)
\\
&\leq C\alpha^{l} \mu(B_{i}).
\end{align*}
We used that $\overline{B_{i}} \cap F\neq \emptyset$, Jensen's inequality and (\ref{t8}), noting that $B_{i}$ is of type 2.\\
If $\underline{B_{i}}$ is of type 1, 
$$
\int_{B_i}   R_{i}^{-l}| b_i|^ld\mu \leq  \int_{B_i}   R_{i}^{-l}| f |^l \leq C\int_{B_{i}} (|\nabla f|^l + |V^{\frac{1}{2}}f|^l)d\mu.
$$
As the same integral but on $\overline{B_{i}}$ is controlled by $\alpha^l \mu( \overline{B_i})$ we get 
$\int_{B_i}   R_{i}^{-l}| b_i|^ld\mu \leq C\alpha^l \mu(B_{i})$.  Since $\nabla b_i = \chi_i \nabla f  + f\nabla\chi_{i}$ we obtain the same bound for
 $\int_{B_i} |\nabla b_i|^ld\mu$. 
Noting that $\overline{B_{i}}\cap F\neq \emptyset$ and  $B_{i}$ is of type 1, we easily deduce that $\int_{B_{i}}|V^{\frac{1}{2}}b_{i}|^{l}\leq C\alpha^{l}\mu(B_{i})$.

 Set $g=f-\sum b_{i}$ where the sum is over both types of balls and is locally finite by 
 \eqref{eqcsds5}. It is clear that
 $g=f$ on $F=M\setminus \Omega$ and $g=\sum{}^2\ f_{B_{i}}  \chi_{i}$ on $\Omega$, where $\sum{}^j$ means that we are summing over cubes of type $j$. Let us prove \eqref{eqcsds1}.
 
First, by the differentiation theorem, $V^{\frac{1}{2}}|f| \leq \alpha$ almost everywhere on $F$.  Next, since
$V \in A_{\infty}$ implies $V^{\frac{l}{2}} \in RH_{\frac{2}{l}}$ we have 
$V_{B_{i}}  \leq C ((V^{\frac{l}{2}})_{B_{i}})^{\frac{2}{l}}$. 
Therefore
$$
\int_{\Omega}  V|g|^2 d\mu\leq \sum{}^2 \int_{B_{i}} V  |f_{B_{i}} |^2
\leq  C \sum{}^2   \left((V^{\frac{l}{2}})_{B_{i}}  ) |f_{B_{i}}|^l \right)^{\frac{2}{l}} \mu(B_i).
$$
Now, by construction of the type 2 balls and the $L^l$ version of Fefferman-Phong inequality,   
$$
(V^{\frac{l}{2}})_{B_{i}} |f_{B_{i}} |^l \leq C( |\nabla  f|^l +  | V^{\frac{1}{2}}f |^l)_{B_{i}} \leq C\alpha^l.
$$
It comes that 
$$
 \int_{\Omega}  V|g|^2d\mu  \leq C \sum{}^2\ \alpha^{2-l}     \mu(B_{i}) \leq C'\alpha^{2-l} \int_{M}(|\nabla f|^l +  | V^{\frac{1}{2}}f |^l)d\mu. 
 $$
 Combining the estimates on $F$ and $\Omega$, we obtain the desired bound for $ \int_{M}  V|g|^2d\mu$. 
 We finish the proof by  estimating  $ \|\,|\nabla g|\,\|_{\infty}$ and $\|\,|\nabla g|\,\|_{l}$. Observe that $g$ is a locally integrable function on $M$. Indeed, let $\varphi\in L_{\infty}$ with compact support. Since $d(x,F)\geq R_{i}$ for $x \in \supp\, \,b_{i}$, we obtain
\begin{equation*} \int\sum_{i}|b_{i}|\,|\varphi|\,d\mu \leq
\Bigl(\int\sum_{i}\frac{|b_{i}|}{R_{i}}\,d\mu\Bigr)\,\sup_{x\in
M}\Bigl(d(x,F)|\varphi(x)|\Bigr).\quad
\end{equation*}
If $B_i$ is of type 2
\begin {align*}
\int \frac{|b_{i}|}{R_{i}}d\mu&\leq \mu(B_i)^{1-\frac{1}{l}}\int \frac{|b_{i}|^{l}}{R_{i}^{l}}d\mu
\\
&\leq C\mu(B_i)^{1-\frac{1}{l}}\int_{B_{i}}|\nabla f|^{l} d\mu
\\
&\leq C\alpha\mu(B_{i}).
\end{align*}
We used the H\"{o}lder inequality, $(P_{l})$ and that $\overline{B_{i}}\cap F\neq \emptyset$, $q'$ being the conjugate of $q$.\\

If $B_i$ is of type 1, 
\begin{equation*}\int \frac{|b_{i}|}{R_{i}}d\mu\leq \mu(B_i)^{1-\frac{1}{l}}\int \frac{|b_{i}|^{l}}{R_{i}^{l}}d\mu\leq C\alpha \mu(B_i).
\end{equation*}
 Hence 
$ \displaystyle \int\sum_{i}|b_{i}||\varphi|d\mu \leq
C\alpha\mu(\Omega)^{\frac{1}{l}} \sup_{x\in M
}\Bigl(d(x,F)|\varphi(x)|\Bigr)$. Since $f\in L^{1}_{loc}$, we conclude that $g\in L^{1}_{loc}$. Thus $\nabla g= \nabla f - \sum \nabla b_{i}$. 
It follows from the $L^{l}$ estimates on $\nabla b_{i}$ and the bounded overlap property that
$$
\left\|\sum |\nabla b_{i}|\right\|_{l} \leq C'(\|\,|\nabla f|\,\|_{l} + \|V^{\frac{1}{2}}f\|_{l}).$$
As $g=f-\sum b_{i}$, the same estimate holds for $\|\,|\nabla g|\,\|_{l}$.
  Next, a computation of the sum $\sum \nabla b_i$ leads us to 
 $$
 \nabla g = \ind _{F} (\nabla f) -\sum{}^1 f\nabla \chi_i- \sum{}^2\ (f-f_{B_{i}} )\ \nabla \chi_{i}.
 $$
 Set $h_{i}=\sum{}^i \ (f-f_{B_{i}})\ \nabla\chi_{i}$ and $h=h_1+h_2$. Then
 $$\nabla g=(\nabla f)\ind_{F}-\sum{}^1 f\nabla \chi_i-(h-h_1)= (\nabla f)\ind_{F}+\sum{}^1f_{B_i}\nabla \chi_i-h.$$
By definition of $F$ and the differentiation theorem, $|\nabla g|$ is bounded by $\alpha$ almost everywhere on $F$. By already seen arguments  for type 1 balls,  $|f_{B_{i}} | \leq C\alpha R_{i} $. Therefore,
$|\sum{}^1f_{B_i}\nabla \chi_i| \leq C\sum{}^1\ \ind_{B_{i}} \alpha\leq CN\alpha$. It remains to control $\|h\|_{\infty}$. For this, note first that $h$ vanishes on $F$ and is locally finite on $\Omega$.
Then fix  $x\in \Omega$. Observe that $\sum_i\nabla \chi_i(x)=0$ and by definition of $I_x$, the sum reduces $i\in I_x$. For all $i\in I_{x}$, we have $|f(x)-f_{B_{i}}|\leq Cr_{i} \alpha$. Hence,  we have for all $j\in I_{x}$,  
$$
\sum_{i}(f(x)-f_{B_{i}})\nabla \chi_i(x) = \sum_{i\in I_{x}}(f(x)-f_{B_{i}})\nabla \chi_i(x)=\sum_{i \in I_{x}}(f_{B_{j}}-f_{B_{i}})\nabla \chi_i(x).
$$
We claim that  $|f_{B_{j}}-f_{B_{i}}|\leq Cr_{j} \alpha$ with $C$ independent of $i,j\in I_{x}$ and $x\in \Omega$.
Indeed, we use that $B_{i}$ and $B_{j}$ are contained in $7 B_{j}$, Poincar\'e inequality $(P_{l})$, the comparability of $r_{i}$ and $r_{j}$, and  that $\overline{B_i}\cap F\neq \emptyset$. Since $I_{x}$ has cardinal bounded by $N$, we are done. 

We conclude that $\|h\|_{\infty}\leq C\alpha$ and interpolating $\|\,|\nabla g|\,\|_{l}$ and $\|\,|\nabla g|\,\|_{\infty}$, we finish therefore the proof. 
\end{proof} 
\begin{prop}\label{DLV} Let $M$ be a complete Riemannian manifold satisfying $(D)$. Let $V\in A_{\infty}$. Moreover assume that $M$ admits a Poincar\'{e} inequality $(P_{p})$ for some $1<p<2$. Then, $\Lip(M)\cap\dot{W}_{2,V^{\frac{1}{2}}}^{1}\cap \dot{W}_{p,V^{\frac{1}{2}}}^{1}$\footnote{$Lip(M)$ is the set of all Lipschitz functions on $M$.} is dense in $\dot{W}_{p,V^{\frac{1}{2}}}^{1}$.
\end{prop}
\begin{proof} Theorem \ref{kz} proves that $M$ admits a Poincar\'{e} inequality $(P_{l})$ for some $1\leq l<p$. Let $f\in \dot{W}_{p,V^{\frac{1}{2}}}^{1}$. For every $n\in \mathbb{N^{*}}$, consider the Calder\'{o}n-Zygmund decomposition of Proposition \ref{CZ} with $\alpha=n$. Take a compact $K$ of $M$. We have
\begin{align*}
\int_{K}|f-g_{n}|^{l}d\mu &= \int_{K\cap (\bigcup_{i}B_{i})} |\sum_{i}b_{i}|^{l}d\mu
\\
&=\int_{\bigcup_{i}K\cap {B_{i}}} |\sum_{i}b_{i}|^{l}d\mu
\\
&\leq C  \sum{}^2\ \int_{K\cap B_{i}}\frac{|f-f_{B_{i}}|^{l}}{R_{i}^{l}} d(x,F_{n})^{l}d\mu+C\sum{}^1\ \int_{K\cap B_{i}}\frac{|f|^{l}}{R_{i}^{l}} d(x,F_{n})^{l}d\mu
\\
& \leq C \sup\limits_{x\in K}(d(x,F_{n}))^{l} \sum_{i}\int_{ B_{i}}(|\nabla f|^{l}+|V^{\frac{1}{2}}f|^{l}) d\mu
\\
& \leq C \sup\limits_{x\in K}(d(x,F_{1}))^{l}\sum_{i}  n^{l} \mu(B_{i})
\\
&\leq C n^{l-p}(\|\,|\nabla f|\,\|_{p}^{p}+\|\,|V^{\frac{1}{2}} f|\,\|_{p}^{p}).
\end{align*}
Letting $n\rightarrow \infty$, we get that $\int_{K}|f-g_{n}|^{l}d\mu \rightarrow 0$. Hence $(f-g_{n})$ converges to $0$ when $n\rightarrow \infty$ in the distributional sense.\\
Let us check that $(V^{\frac{1}{2}}(f-g_{n}))_n$ is bounded in $L^{p}$. Indeed,
\begin{align*}
\int_{M}|V^{\frac{1}{2}}(f-g_{n})|^{p}d\mu &\leq\int_{\Omega_{n}}|V^{\frac{1}{2}}f|^{p}d\mu+\sum{}^{2}\int_{\Omega_{n}} V^{\frac{p}{2}}|f_{B_{i}}|^{p}d\mu
\\
&\leq\int_{\Omega_{n}}|V^{\frac{1}{2}}f|^{p}d\mu+\sum{}^{2} \left((V^{\frac{l}{2}})_{B_{i}}|f_{B_{i}}|^{l}\right)^{\frac{p}{l}}\mu(B_{i})
\\
&\leq \int_{\Omega_{n}}|V^{\frac{1}{2}}f|^{p}d\mu+Cn^{p}\mu(\Omega_{n})
\\
&\leq C(\|\,|\nabla f|\,\|_{p}^{p}+\|V^{\frac{1}{2}}f\|_{p}^{p}).
\end{align*}
 Similarly,
$$
\int_{M}|\nabla f-\nabla g_{n}|^{p}d\mu =\int_{\Omega_{n}}|\nabla f-\nabla g_{n}|^{p}d\mu \leq C\int_{\Omega_{n}}|\nabla f|^{p} d\mu+ Cn^{p}\mu(\Omega_{n})\leq C.
$$
Thus, $(\nabla f-\nabla g_{n})_{n}$ is bounded in $L^{p}$. So $(f-g_{n})_{n}$ is bounded in $\dot{W}_{p,V^{\frac{1}{2}}}^{1}$. Since $\dot{W}_{p,V^{\frac{1}{2}}}^{1}$ is reflexive --Proposition \ref{RW}--, there exists a subsequence, which we denote also by $(f- g_{n})_{n}$, converging weakly in $\dot{W}_{p,V^{\frac{1}{2}}}^{1}$ to a function $h$. 
The uniqueness of the limit in the distributional sense yields $h=0$. By Mazur's Lemma, we find a sequence $(h_{n})$ of convex combinations of $(f-g_{n})$ such that $h_{n}=\sum_{k=1}^{n} a_{n,k}(f-g_{k})$, $a_{n,k}\geq0$, $\sum_{k=1}^{n}a_{n,k}=1$, that converges to $0$ in $\dot{W}_{p,V^{\frac{1}{2}}}^{1}$. Since $\nabla h_{n}=\nabla f-\nabla l_{n}$ and $V^{\frac{1}{2}} h_{n}=V^{\frac{1}{2}} (f- l_{n})$ with $l_{n}=\sum_{k=1}^{n} a_{n,k}g_{k}$, we obtain $l_{n}\underset {n\rightarrow \infty}{\longrightarrow}f$ in $\dot{W}_{p,V^{\frac{1}{2}}}^1$ and the proposition follows on noting that $g_{n}$, hence $l_{n}$,  also belongs to $Lip(M)\cap \dot{W}_{2,V^{\frac{1}{2}}}^{1}$.
\end{proof}

\subsection{Estimates for subharmonic functions}
Fix an open set $\Omega\subset M$. A subharmonic function on $\Omega$ is a function $v\in L^{1}_{loc}(\Omega)$ such that $\Delta v \geq 0$ in $D'(\Omega)$. 
\begin{lem}\label{Eu}  Let $M$ be a Riemannian manifold satisfying $(D)$ and $(P_2)$. Let $R>0$ and $x_{0}$ be a point such that a neighborhood of $\overline{B(x_{0},4R)}$ is contained in $M$.
Suppose that $f$ is a non-negative subharmonic function  defined on this neighborhood. Then, there is a constant $C>0$ independent of $f$, $x_0$, $R$ such that 
\begin{equation}\label{u}
\sup\limits_{x\in B(x_{0},R)}f(x)\leq C\left(\aver{B(x_{0},4R)}f^2(y)d\mu(y)\right)^{\frac{1}{2}}
\end{equation}
It readily follows from Lemma \ref{buck2} that for all $r>0$, $1<\eta< 4$, there is $C>0$ such that
\begin{equation}\label{ug}
 \sup\limits_{x\in B(x_{0},R)}f(x)\leq C\left(\aver{B(x_{0},\eta R)}f^{r}(y)d\mu(y)\right)^{\frac{1}{r}}.
\end{equation}
\end{lem}
\begin{proof}
In \cite{lip2}, Theorem 7.1, this lemma is stated for Riemannian manifolds with non-negative Ricci curvature. The proof relies on the following properties of the manifold. First, the Harnack inequality for non-negative harmonic functions which holds for complete Riemannian manifolds satisfying $(D)$ and  $(P_2)$ (see \cite{grigoryan1}). Secondly, the Poincar\'e inequality $(P_2)$. Finally, the Caccioppoli inequality for non-negative subharmonic functions --Lemma 7.1 in \cite{lip2}-- which is valid on any complete Riemannian manifold. We then get this lemma under the hypotheses $(D)$ and $(P_2)$.
\end{proof}
Other forms of the mean value inequality for subharmonic functions still hold if the volume form is replaced 
by a weighted measure of Muckenhoupt type. More precisely, 

\begin{lem}\label{subw} Consider a complete Riemannian manifold $M$ satisfying $(D)$ and $(P_{2})$. Let $V\in A_{\infty}$ and $f$ a non-negative subharmonic function defined on a neighborhood of $\overline{B(x_{0},4R)}$, $0<s<\infty$ and $1<\eta< 4$. Then for some $C$ depending on the $A_\infty$ constant of $V$, $s$ (and independent of $f$ and $x_0,\,R$),  we have
$$
\sup_{x\in B(x_{0},R)}f(x) \leq \left(\frac {C}{V(B(x_{0},\eta R))} \int_{B(x_{0},\eta R)} V f^sd\mu\right)^{\frac{1}{s}}.$$
\end{lem}

Here $V(E)=\int_{E}Vd\mu$. 
As $A_{\infty}$ weights have the doubling property we have $V_{B(x_{0},\eta R)}\sim V_{B(x_{0},R)}$ and the inequality above
is the same as
\begin{equation}
\label{eq:subbis}
V_{B(x_{0},R)} (\sup_{B(x_{0},R)}f^s) \leq C ( Vf^s)_{B(x_{0},\eta R)}.
\end{equation}

\begin{proof} Since $V\in A_{\infty}$, there is $t<\infty$ such that  $V\in A_{t}$. Hence
for any  non-negative measurable function $g$ we have 
\begin{align*}
g_{B(x_{0},\eta R)} &\leq C \left(\frac{1}{V(B(x_{0},\eta R))} \int_{B(x_{0},\eta R)} V  g^td\mu \right)^{\frac{1}{t}}
\\
&= C\left((V g^t )_{B(x_{0},\eta R)}\right)^{\frac{1}{t}} \left( V_{B(x_{0},\eta R)}\right)^{-\frac{1}{t}} .  
\end{align*}
 Applying (\ref{ug})  with $r=\frac{s}{t}$   yields
$$
f(x)\leq C \left((f^{\frac{s}{t}})_{B(x_{0},\eta R)}\right)^{\frac{t}{s}} \leq  C  \left(  (Vf^s )_{B(x_{0},\eta R)}\right)^{\frac{1}{s}} \left(  V_{B(x_{0},\eta R)}\right)^{-\frac{1}{s}} .
$$
\end{proof}
\begin{cor}\label{lem:rh2q} 
Let $M$ be a complete Riemannian manifold satisfying $(D)$ and $(P_{2})$. Let $V\in RH_{r}$ for some $1<r\leq \infty$, $0<s<\infty$ and $1< \eta\leq4$. Then there is  $C\geq 0$  depending only on the $RH_r$ constant of $V$, $s$  such that  for any  ball $B(x_{0},R)$  and any non-negative subharmonic function defined on a neighborhood of $\overline{B(x_{0},4R)}$  we have 
$$
\left( ((Vf^s)^{r})_{B(x_{0},R)} \right)^{\frac{1}{r}} \leq C ( Vf^s)_{B(x_{0},\eta R)}.
$$
\end{cor}

\begin{proof}  We have

$$
\left( ((Vf^s)^{r})_{B(x_{0},R)} \right)^{\frac{1}{r}} \leq C  \left( (V^r)_{B(x_{0},R)} \right)^{\frac{1}{r}} \sup_{B(x_{0},R)} f^s \leq 
C V_{B(x_{0},R)} \sup_{B(x_{0},R)} f^s  \leq C (Vf^s)_{B(x_{0},\eta R)}.
$$
The second inequality uses the $RH_{r}$ condition on $V$ and the last inequality is (\ref{eq:subbis}).
\end{proof}
\
\section{Maximal inequalities}
This section is devoted to the proof of Theorem \ref{th:1}. Let $1<q\leq\infty$ and $V\in RH_{q}$. 
The following lemma is classical in an Euclidean setting \cite{gallouet}, \cite{kato} (see also \cite{auscher3}).

\begin{lem}\label{lem:l1} Let $M$ be a complete Riemannian manifold. We assume that $V\in L_{loc}^{1}(M)$ is not identically 0. 
Let $u\in C^{\infty}_{0}(M)$. Then
 $$
\int_{M} V|u| d\mu \leq \int_{M} |(-\Delta+V)u|d\mu, 
$$
 $$
\int_{M} |\Delta u| d\mu \leq 2 \int_{M}|(-\Delta+V)u|d\mu. 
$$
\end{lem}
 \begin{proof}
Let us prove the estimate for $V|u|$. 
Take $p_{n}:\mathbb{R}\rightarrow \mathbb{R}$ a sequence of $C^1$ functions such that $|p_{n}|\leq C,\, p'_{n}(t)\geq 0$ and $p_{n}(t)\rightarrow \,sign(t)$ for every $t\in \mathbb{R}$. Using the Lebesgue convergence theorem we see that
$$
-\int_{M}\,sign(u)\Delta u d\mu=-\lim_{n}\int_{M}p_{n}(u)\Delta u d\mu=\lim_{n}\int_{M}|\nabla u|^{2}p'_{n}(u)d\mu\geq 0.
$$
If $-\Delta u+Vu=f$, we get
$$
\int_{M}V|u|d\mu\leq \int_{M}\,sign(u)(-\Delta+V)ud\mu=\int_{M}f\,sign(u)d\mu\leq \int_{M}|f|d\mu.
$$
This gives the desired estimation for $V|u|$.

The estimate for $\Delta u$ follows from that of $Vu$ since $-\Delta u+Vu=f$.
\end{proof}
Let $\mathcal{D}_{1}(H)=\{ u\in L^1_{loc}\, ; \, Vu\in L^1_{loc},\, (-\Delta+V) u \in L^1 \}$.  One can easily check that $C_{0}^{\infty}$ is dense in $\mathcal{D}_{1}(H)$ (\cite{carbonaro} for a proof in the Euclidean parabolic case)  thanks to the Kato inequality on manifolds (\cite{braverman}, Theorem 5.6). 
 Thus the above estimates for $\int V|u|$ and $\int |\Delta u|$ still holds for any $u\in \mathcal{D}_{1}(H)$. Lemma \ref{lem:l1} shows that $\mathcal{D}_{1}(H) =\{ u\in L^1_{loc}\, ; \, \Delta u \in L^1,  Vu\in L^1\}$
equipped with the topology defined by the semi-norms for $L^1_{loc}$, 
$\|\Delta u \|_{1}$ and $\|Vu\|_{1}$. 
 We have therefore obtained
\begin{thm}\label{th:L1} The operator $H^{-1}$ \textit{a priori} defined on  
$L^\infty_{0}(M)$ --the set of compactly supported bounded functions defined on $M$-- extends to a bounded operator from $L^1(M)$ into $\mathcal{D}_{1}(H)$.
Denoting again $H^{-1}$ this extension, $V\, H^{-1}$ is a  positivity-preserving
contraction on $L^{1}(M)$ and $\frac 1 2 \Delta H^{-1}$ is a contraction on
$L^{1}(M)$.
\end{thm}

\begin{prop}\label{pro:uniqueness} Assume that $M$ satisfies $(D)$ and $(P_2)$. Let $f\in L^1(M)$.
There is uniqueness of solutions for the equation $-\Delta u + Vu =f $ in the class $L^1(M)\cap \mathcal{D}_1(H)$. In particular, if $u \in C_{0}^\infty(M)$ and $f=-\Delta u
+ Vu$, then 
$u=H^{-1}f$. 
\end{prop}

\begin{proof} Assume $-\Delta u + Vu =0$, then for $\epsilon>0$  we have $-\Delta
u + Vu + \epsilon u = \epsilon u $. As $u\in L^1(M)$, we can write 
$|u| \le (-\Delta +\epsilon)^{-1}(\epsilon | u|)= (-\epsilon^{-1} \Delta +
1)^{-1}| u|$. Using the upper bound of the kernel of $(-\epsilon^{-1}
\Delta + 1)^{-1}$ which follows from $(D)$ and $(P_2)$, and taking limits when $\epsilon \to 0$ we get $u=0$. 
\end{proof}

\begin{cor}\label{cor:L1} Assume $(D)$ and $(P_2)$. Then equation (\ref{max}) holds.
\end{cor}
\begin{proof}   If $u \in C_{0}^\infty(M)$ and $f=-\Delta u + Vu$, then 
$Vu=V H^{-1}f$ and $\Delta u = \Delta  H^{-1} f$ by the proposition above. 
Applying Theorem \ref{th:L1} we get 
$\|Vu \|_{1} \le  \|-\Delta u + Vu \|_{1}$ and $\|\Delta u\|_{1}\le 2 \|-\Delta
u + Vu \|_{1}$.
\end{proof}
 We now give the following criterion for $L^{p}$ boundedness:
\begin{thm} (\cite{auscher7}) \label{theor:shen} Let $M$ be a complete Riemannian manifold satisfying $(D)$. 
Let $1\leq p_{0} <q_0\le \infty$. Suppose that $T$ is a bounded
sublinear operator  on $L^{p_{0}}(M)$. Assume that there exist
constants $\alpha_{2}>\alpha_{1}>1$, $C>0$ such that
\begin{equation}\label{T:shen}
\big(\aver{B} |Tf|^{q_0}\big)^{\frac{1}{q_0}}
\leq
C\, \bigg\{ \big(\aver{\alpha_{1}\, B}
|Tf|^{p_0}\big)^{\frac{1}{p_0}} +
(S|f|)(x)\bigg\},
\end{equation}
for any ball $B$, $x\in B$ and  all $f\in L^{\infty}_{0}(M)$ with 
support in $M \setminus \alpha_{2}\, B$, where $S$ is a positive operator.
Let $p_{0}<p<q_{0}$. If $S$ is bounded on $L^p(M)$, then, there is a constant $C$ such that
$$
\|T f\|_{p}
\leq
C\, \|f\|_{p}
$$
for all $f\in L_{0}^{\infty}(M)$.
\end{thm}

Note that the space $L_{0}^{\infty}(M)$ can be replaced by $C_{0}^\infty(M)$.
\\
Now we use the $L^1$ estimate and Theorem \ref{theor:shen}  to get

\begin{thm}\label{Maxp}
 Let $M$ be a complete Riemannian manifold satisfying $(D)$ and $(P_2)$. Consider
$V\in RH_{q}$, with $q>1$. Then, there exists $r>q$, such
that $VH^{-1}$ and $\Delta H^{-1}$ defined on $L^{1}(M)$ by Theorem \ref{th:L1} extend to $L^p(M)$ bounded operators  for all
$1<p<r$.
\end{thm}

\begin{proof} By difference, it suffices to prove the theorem for $V H^{-1}$. We
know that this is a bounded operator on $L^1(M)$. Let $r$ be given by the self-improvement of the reverse H\"{o}lder condition of $V$. Fix a ball $B$ and let $f\in
L^\infty(M)$ with compact support contained in $M \setminus 4B$. Then $u=H^{-1}
f$ is well-defined in  $\dot{\mathcal{V}}$ and is a weak solution of $-\Delta
u + Vu=0$ in $4B$. Since $|u|^2$ is subharmonic (cf section 8.1), we can apply Corollary \ref{lem:rh2q} with
$V$, $f=|u|^2$ and $s=\frac{1}{2}$. Thus \eqref{T:shen} holds with $T= VH^{-1}$, 
$p_{0}=1$, 
$q_{0}=r $, $S=0$, $\alpha_{1}=2$ and $\alpha_{2}=4$. Hence, Theorem \ref{theor:shen} asserts that $T=VH^{-1}$ is  bounded on
$L^p(M)$ for $1<p<r$. 
\end{proof}

\paragraph{\bf Proof of Theorem \ref{th:1} } 
  Let $u\in C_{0}^{\infty}(M)$ and $f=-\Delta u + Vu$. Proposition \ref{pro:uniqueness} shows that $u=H^{-1}f$. Since  $V\in RH_{q}$, Theorem \ref{Maxp} shows that $VH^{-1}$ and $\Delta H^{-1}$ have bounded extensions on $L^p(M)$ for $1<p<q+\epsilon$ for some $\epsilon>0$ depending on $V$. 
This means that $\|Vu \|_{p} +\|\Delta u\|_{p}\lesssim \|f \|_{p}$ which is the desired result.   
 

\section{Complex interpolation} 

We shall use complex interpolation to obtain item 1. of Theorem \ref{RH}. This method
is based on the boundedness of imaginary powers of $H$ and of the Laplace-Beltrami
operator. Then we use Stein's interpolation theorem to prove the boundedness
of $\nabla H^{-\frac{1}{2}}$ and $V^{\frac{1}{2}} H^{-\frac{1}{2}}$ on $L^p(M)$ for $1<p<2(q+\epsilon)$ and therefore obtain item 1. of Theorem \ref{RH}.
 
 Let $y\in \mathbb{R}$, the operator $H^{iy}$ is defined  via spectral theory.
One has
$$\|H^{iy}\|_{2\rightarrow 2}  =1.$$


\begin{thm}\label{HI} Let $M$ be a complete Riemannian manifold satisfying $(D)$ and assume that the heat kernel verifies the following upper bound: for all $x\in M$ and $t>0$
\begin{equation}\label{Up}
p_{t}(x,x)\leq \frac{C}{\mu(B(x,\sqrt{t}))}.
\end{equation}
Let $V$ be a non-negative locally integrable function on $M$. Then for all $\gamma\in
\mathbb{R}$,  $H^{i\gamma}$  has a bounded extension on $L^p(M)$, $1<p<\infty$, and
for fixed $p$ its operator  norm does not exceed
$C(\delta,p )e^{\delta |\gamma|}$ for some $\delta >0$.  
\end{thm}
\begin{rem}
The operator norm is far from optimal but sufficient for us.
\end{rem}
\begin{proof}[Proof of Theorem \ref{HI}] For $V=0$, this follows from the universal multiplier theorem for Markovien semi groups (Corollary 4, p.121 in \cite{stein4}). However, the following proof works for all $V$. Indeed, the remark after Theorem 3.1 in \cite{duong} applies to $H$: $H$ has a bounded holomorphic functional calculus on $L^{2}(M)$ in any sector $|arg z|<\theta$, $0<\theta<\pi$ and the kernel $h_{t}(x,y)$ of $e^{-tH}$ has a Gaussian upper bound. This follows from the domination of $e^{-tH}$ by $e^{-t\Delta}$, $(D)$ and $(\ref{Up})$. We have
$$
|h_{t}(x,y)|\leq \frac{C}{\mu(B(x,\sqrt{t}))}e^{-c\frac{d^{2}(x,y)}{t}}
$$
for every $t>0$, $x,\,y\in M$.

Thus  a variant of Theorem 3.1 in \cite{duong} (see page 104 there) shows that $H$ has a bounded holomorphic functional calculus on $L^{p}(M)$ in any sector $|arg z|<\mu, \,\frac{\pi}{2}<\mu\leq \pi$ for $1<p<\infty$. This implies
$$
\|H^{i\gamma}\|_{p\rightarrow p}\leq C(p,\mu) \sup_{|argz|<\mu}|z^{i\gamma}|\leq C_{p,\mu}e^{|\gamma|\mu}.
$$   
\end{proof}
\begin{lem}
The space $\mathcal{D}=\mathcal{R}(H) \cap L^1(M)\cap L^\infty(M)$ is dense in
$L^p(M)$ for $1<p<\infty$. 
\end{lem}
\begin{proof} Same proof as that of Lemma 6.2 in \cite{auscher3}.
\end{proof}
\begin{prop} \label{alpha} Assume that $M$ satisfies $(D)$ and $(P_2)$. Let $V\in RH_{q}$ for some $1\leq q<\infty$. Then, for $0<\alpha<1$, there exists $\epsilon>0$ such that the operators $\Delta^{\alpha}H^{-\alpha}$, $V^{\alpha}H^{-\alpha}$ are bounded on $L^{p}(M)$ for $1<p<\frac{1}{\alpha}q+\epsilon$.
\end{prop}

\begin{proof} From Theorem \ref{HI}, we have that $\Delta^{i\gamma}$ and $H^{i\gamma}$ are $L^{p}(M)$ bounded for  $1<p<\infty$ and  $\gamma \in \mathbb{R}$. Moreover, Theorem \ref{th:1} asserts that $\Delta H^{-1}$ and $VH^{-1}$ are $L^{p}(M)$ bounded for $1<p<q+\epsilon$ for some $\epsilon >0$. It follows from Stein's interpolation theorem  \cite{stein3} that $\Delta^{\alpha}H^{-\alpha}$, $V^{\alpha}H^{-\alpha}$ are bounded on $L^{p}(M)$ for $1<p<\frac{1}{\alpha}(q+\epsilon)$ (see \cite{auscher3} for details).
\end{proof}

We can now prove item 1. of Theorem \ref{RH}. Fix $1<p<2(q+\epsilon)$. Let $u\in C_{0}^\infty(M)$. 
Since $u \in
\mathcal{V}$, $f=H^{\frac{1}{2}}u$ is well-defined. We assume that $f\in L^p(M)$, 
otherwise there is nothing to prove. Applying Proposition \ref{alpha} to $V^{\frac{1}{2}}$, it comes that $\|V^{\frac{1}{2}}u\|_{p}\leq C_p \|f\|_p$. 
 The $L^{p}(M)$ boundedness of the Riesz
transform which holds for all $1<p<p_{0}$ with $p_0>2$ on a complete Riemannian manifold satisfying $(D)$ and $(P_2)$ and again Proposition \ref{alpha} yield
$$
\|\,|\nabla u|\, \|_{p}+ \le C(p) \|\Delta^{\frac{1}{2}} H^{-\frac{1}{2}} f\|_{p} \le C'(p) \|f\|_{p}$$
for $1<p<\inf(p_0,2(q+\epsilon))$
and finishes the proof. 
\begin{rem}
This interpolation argument also gives us a proof of the $L^p(M)$ boundedness of 
$\nabla H^{-1}$ and $V^{\frac{1}{2}}H^{-\frac{1}{2}}$ for $1<p<2$ for all non zero $V\in
L^1_{loc}(M)$.
\end{rem}
  \section{Proof of Theorem \ref{RRT}}
 The proof is similar to that of item 2. of Theorem 1.2 in \cite{auscher3} with some modifications. We write it for the sake of completeness.
Denote $H=-\Delta+V$. Assume that $1<l<2$. Let $f\in \Lip(M)\cap \dot{W}_{l,V^{\frac{1}{2}}}^{1}\cap\dot{W}_{2,V^{\frac{1}{2}}}^{1}$.
 We use the following resolution of $H^{\frac{1}{2}}$:
$$
H^{\frac{1}{2}}f= c\int_0^\infty H e^{-t^2 H} f \, dt
$$
in the distributional sense.
It suffices to obtain the 
result for the truncated integrals $\int_\epsilon^R\ldots$ with bounds independent of $\epsilon,R$, 
and then to let
$\epsilon\searrow 0$ and $R\nearrow \infty$. 
For the truncated integrals, all the calculations are justified. We thus consider that $H^{\frac{1}{2}}$ is one of the truncated integrals but we still write the limits as $0$ and $+\infty$ to simplify the exposition. As $f$ does not belong to $C_{0}^{\infty}(M)$, we have to give a meaning to $He^{-tH}f$ for $t>0$. 
 Take $\eta_{r}$ a smooth function on $M$, $0\leq \eta_r\leq 1$, $\eta =1$ on a ball $B$ of radius $r>0$, $\eta_r=0$ outside $2B$ and $\|\,|\nabla \eta_r|\,\|_{\infty}\leq \frac{C}{r}$. For $\varphi \in C_{0}^{\infty}(M)$,
 \begin{align*}
 \int_{M} f\,H e^{-t^{2}H}\varphi d\mu=&=\lim_{r\rightarrow \infty}\int_{M}\eta_{r}fH e^{-t^{2}H}\varphi d\mu
 \\
 &=\int_{M}\eta_{r}\nabla f. \nabla e^{-t^{2}H}\varphi d\mu+\int_{M}f\nabla \eta_{r}.\nabla e^{-t^{2}H}\varphi d\mu
 \\
 &+\int_{M}\eta_{r}f \, V e^{-t^{2}H}\varphi d\mu
 \\
 &=I_r+II_r+III_r.
 \end{align*}
 We used Fubini and Stokes theorems.
 Note that $\int_{M}|\nabla_{x}h_{t}(x,y)|^{2}e^{\gamma \frac{d^{2}(x,y)}{t}}d\mu(x)\leq \frac{C}{t\mu(B(y,\sqrt{t}))}$. This is due to the Gaussian upper estimate of the kernel $h_{t}$ of $e^{-tH}$ and that of $\partial_{t}h_{t}$ under $(D)$ and $(P_2)$ (see \cite{coulhon2}, Lemma 2.3, for the heat kernel $p_t$ of $e^{-t\Delta}$).
Since $|\nabla f|\in L^2(M)$ then $I_r\rightarrow\int_{M}\nabla f.\nabla e^{-t^{2}H}\varphi d\mu$. Since $f$ is Lipschitz, $II_r\rightarrow 0$.
We have also $\int_{M}|h_t(x,y)|^2e^{\gamma\frac{d^{2}(x,y)}{t}}d\mu(x)\leq\frac{C}{\mu(B(y,\sqrt{t}))}$ and  $V^{\frac{1}{2}}f\in L^2(M)$. Thus $III_r \rightarrow \int_{M}fVe^{-t^{2}H}\varphi d\mu$.  
 This proves that $He^{-t^{2}H}f$ is defined as a distribution by
 $$\langle He^{-t^2H}f,\varphi\rangle=\int_{M}\nabla f.\nabla e^{-t^2H}\varphi d\mu+\int_{M}V^{\frac{1}{2}}fV^{\frac{1}{2}}e^{-t^2H}\varphi d\mu.$$
 Therefore, integrating in $t$ yields
 $$\langle H^{\frac{1}{2}}f,\varphi\rangle=\langle\nabla f, \nabla H^{-\frac{1}{2}}\varphi\rangle+\langle V^{\frac{1}{2}}f,H^{-\frac{1}{2}}\varphi\rangle.
 $$
 
\ 
 We return to the proof of Theorem \ref{RRT}. Apply the Cal\-de\-r\'on-Zygmund decomposition of Lemma \ref{CZ} 
to $f$  at height $\alpha$ and write $f=g+\sum_i b_i$.

For $g$, we have
 \begin{align*}\label{eq22}
\mu\left(\left\lbrace x \in M; |H^{\frac{1}{2}}g(x)| >\frac {\alpha}{ 3}\right\rbrace \right)&\leq 
\frac{9}{\alpha^2}\int
|H^{\frac{1}{2}} g|^2d\mu \leq \frac{9}{\alpha^2}\int(|\nabla g|^2 + V|g|^2)d\mu 
\\
&\leq 
\frac{C}{\alpha^l}\int (|\nabla f|^l +  |V^{\frac{1}{2}}f|^l)d\mu.
\end{align*}
We used a similar argument as above to compute $H^{\frac{1}{2}}g$ (see \cite{auscher1}) and the $L^2$ estimate follows.
 For the last inequality we used \eqref{eqcsds2} of the Calder\'on-Zygmund decomposition and that $l<2$.
 
The argument to estimate  $H^{\frac{1}{2}}b_i$ will use the Gaussian upper bound of $h_{t}$. As we mentioned above, under our assumptions we have the Gaussian upper bound for the kernel of $e^{-t^{2}H}$ and by analyticity for $He^{-t^{2}H}$. As $b_{i}$ is supported in a ball and integrable $He^{-t^{2}H}b_i$ is defined by the convergent integral $\int_{M}\frac{-1}{2t}\partial_{t}h_{t^2}(x,y)b_i(y)d\mu(y)$.
 Let $r_i=2^k$ if $2^k \leq R_i < 2^{k+1}$ ($R_{i}$ is the radius of $B_{i}$) and
set 
$T_i= \int_0^{r_i} He^{-t^{2} H}  \, dt$ and
 $U_i= \int_{r_i}^\infty He^{-t^2 H}  \, dt$. It is enough to 
estimate $$A=\mu\left(\left\lbrace x \in M; |\sum_i T_ib_i(x)| >\frac{\alpha}{3}\right\rbrace\right)$$ and 
$$
B=\mu\left(\left\lbrace x \in M; \bigg|\sum_i U_ib_i(x)\bigg| >\frac{\alpha}{3}\right\rbrace\right).$$

First,  
$$A \leq \mu(\bigcup_i \overline{B_i}) + \mu\left(\left\lbrace x \in M \setminus \bigcup_i  \overline{B_i} ; 
\bigg|\sum_i T_ib_i(x)\bigg| >\frac {\alpha}{ 3}\right\rbrace\right),$$
and by \eqref{eqcsds4}, $\mu(\bigcup_i \overline{B_i}) \leq \frac{C}{\alpha^l}\int (|\nabla f|^l+  |V^{\frac{1}{2}}f|^l)d\mu$.

For the other term, we have
$$
\mu\left(\left\lbrace x \in M \setminus \bigcup_i \overline{B_i} ; \bigg|\sum_i
T_ib_i(x)\bigg| >\frac {\alpha}{ 3}\right\rbrace\right)
\leq
\frac{C}{\alpha^2}\int \bigg|\sum_i h_i\bigg|^2$$
 with $h_i = \ind_{(\overline{B_i})^c}|T_ib_i|$. To estimate the $L^2$ norm, we dualize against
$u\in L^{2}(M)$ with $\|u\|_{2}=1$:
$$
\int| u| \sum_i h_i = \sum_i\sum_{j=2}^\infty A_{ij}
$$  
where 
$$
A_{ij}= \int_{C_{j}(B_i)} |T_ib_i||u|d\mu, \quad C_{j}(B_{i}) = 2^{j+1}B_{i} \setminus 2^jB_{i}.
$$   
By Minkowski integral inequality, for some appropriate  positive constants $C,c$,
$$
\| T_ib_i \|_{L^2(C_j(B_i))} \leq   \int_0^{r_i} \|He^{-t^{2}H}  b_i \|_{L^2(C_j(B_i))} \, dt.
$$
By the well-known Gaussian upper bounds for the kernels of  
 $tHe^{-tH}$, $t>0$, valid since we have $(D)$ and $(P_{2})$
 $$
 |H e^{-t^{2}H}b_{i}(x)|\leq \int_{M}\frac{C}{t^{2}\mu(B(y,t))}e^{-\frac{cd^{2}(x,y)}{t^{2}}}|b_{i}(y)|d\mu(y).
 $$
 Now $y\in \supp\,b_{i}$, that is $B_{i}$, and $x\in C_{j}(B_{i})$, hence one may replace $d(x,y)$ by $2^{j}r_{i}$ in the Gaussian term since $r_{i}\sim R_{i}$. Also if $x_{i}$ denotes the center of $B_{i}$, we have
 $$
 \frac{\mu(B(x_{i},t))}{\mu(B(y,t))}= \frac{\mu(B(x_{i},t))}{\mu(B(x_{i},r_{i}))} \frac{\mu(B(x_{i},r_{i}))}{\mu(B(y,r_{i}))} \frac{\mu(B(y,r_{i}))}{\mu(B(y,t))}.
 $$
 By $(D)$ and Lemma \ref{CD} 
 as $t\leq r_{i}$, we have
 $$
 \frac{\mu(B(x_{i},t))}{\mu(B(y,t))}\leq C(2\frac{r_{i}}{t})^{s}. 
$$
Using the estimate \eqref{eqcsds3}, 
$\|b_i\|_1\leq
 c
\alpha
R_i \mu(B_{i})$, and $\mu(B_{i})\sim\mu(B(x_{i},r_{i}))$, it comes that
\begin{align*}
|He^{-t^{2}H}b_{i}(x)|&\leq \frac{C}{t^{2}\mu(B(x_{i},t))}\left(\frac{r_{i}}{t}\right)^{s}e^{-\frac{c4^{j}r_{i}^{2}}{t^{2}}}\int_{B_{i}}|b_{i}|d\mu
\\
&\leq \frac{Cr_{i}}{t^{2}}\left(\frac{r_{i}}{t}\right)^{2s}e^{-\frac{c4^{j}r_{i}^{2}}{t^{2}}}\alpha.
\end{align*}
Thus 
$$\|H e^{-t^{2} H}  b_i \|_{L^2(C_j(B_i))} \leq \frac{Cr_{i}}{t^{2}}\left(\frac{r_{i}}{t}\right)^{2s} e^{-\frac{c 4^jr_i^2}{t^2}}\mu(2^{j+1}B_{i})^{\frac{1}{2}}\alpha.
$$
Plugging this estimate inside the integral, we get
$$
\| T_ib_i \|_{L^2(C_j(B_i))} \leq   C\alpha  e^{-c4^j}  \mu(2^{j+1}B_{i})^{\frac{1}{  2}}.
$$
Now remark that for any $y \in B_i$ and any $j\geq 2$, 
$$
\left( \int_{C_j(B_i)} |u|^{2}\right)^{\frac{1}{2}} \leq \left( \int_{2^{j+1}B_i} |u|^{2}\right)^{\frac{1}{2}} \leq  (2^{s(j+1)}\mu(B_i))^{\frac{1}{2}}
\big(\mathcal{M}(|u|^{2})(y)\big)^{\frac{1}{2}}.
$$
Applying H\"older inequality, one obtains
$$
A_{ij} \leq  C\alpha 2^{sj} e^{-c4^j}  \mu(B_{i}) \big(\mathcal{M}(|u|^{2})(y)\big)^{\frac{1}{2}}.
$$
Averaging over $B_i$ yields
$$
A_{ij} \leq C\alpha 2^{sj} e^{-c4^j} \int_{B_i} \big(\mathcal{M}(|u|^{2})(y)\big)^{\frac{1}{2}}\, d\mu(y).$$
Summing over $j\geq 2$ and $i$, it follows that
$$
\int| u| \sum_i h_i d\mu \leq C \alpha \int \sum_i \ind_{B_i}(y) \big(\mathcal{M}(|u|^{2})(y)\big)^{\frac{1}{2}}\, d\mu(y).$$
Using finite overlap \eqref{eqcsds5} of the balls $B_i$ and Kolmogorov's inequality, one obtains
$$
\int| u| \sum_i h_i d\mu \leq C'N\alpha \mu( \bigcup_i B_i )^{\frac{1}{2}} \||u|^2\|_1^{\frac{1}{2}}.
$$
Hence 
$$
\mu\left(\left\lbrace x \in M \setminus \bigcup_i \overline{B_i} ; \bigg|\sum_i T_ib_i(x)\bigg| >\frac \alpha 3\right\rbrace\right) \leq C \mu( \bigcup_i B_i) \leq 
\frac{C}{\alpha^{l}}\int (|\nabla f|^{l}+ |V^{\frac{1}{2}}f|^{l})d\mu$$
by \eqref{eqcsds5} and \eqref{eqcsds4}.

It remains to handle the term $B$. Using functional calculus for $H$ one can
compute 
$U_i$ as $r_i^{-1}\psi(r_i^2 H)$ with $\psi$ the holomorphic function on the sector $|\arg z \,| < {\frac \pi 2}$ given
by 
\begin{equation*}\label{eqpsi}
\psi(z)= \int_1^\infty e^{-t^2z} z\, dt.
\end{equation*}
It is easy to show that $|\psi(z)| \leq C|z|^{\frac{1}{2}} e^{-c|z|}$, uniformly on subsectors $|\arg z\,  | \leq \mu < {\frac \pi 2}$.

The $(P_{l})$ Poincar\'e inequality gives us if $B_{i}$ is of type 2 

$$\|b_i\|_l^{l} \leq C R_i^{l}\int_{B_{i}}|\nabla f|^{l}d\mu\leq C R_{i}^{l}\alpha^{l}\mu(B_{i}).
$$
If $B_{i}$ is of type 1  
\begin{equation}\label{bi}
 b_{i}= (b_{i}-(b_{i})_{B_{i}})\ind_{B_{i}}+(b_{i})_{B_{i}}\, \ind_{B_{i}}.
\end{equation}
Therefore using the type 1 property of $B_{i}$ and also (\ref{bi}) yield
\begin{align*}
\int_{B_{i}}|b_{i}|^{l}d\mu&\leq 2^{l-1}\left(\int_{B_{i}}|b_{i}-(b_{i})_{B_{i}}|^{l}+\mu(B_{i})\,|\aver{B_{i}}b_{i}d\mu|^{l}\right)
\\
&\leq C R_{i}^{l}\mu(B_{i})^{1-l}\int_{B_{i}}|\nabla b_{i}|^{l}d\mu+ C\mu(B_{i})R_{i}^{l}\aver{B_{i}}(|\nabla f|^{l}+|V^{\frac{1}{2}}f|^{l})d\mu
\\
&\leq CR_{i}^{l}\mu(B_{i})^{1-l}\int_{B_{i}}|\nabla f|^{l}d\mu+ C\mu(B_{i})R_{i}^{l}\int_{B_{i}}(|\nabla f|^{l}+|V^{\frac{1}{2}}f|^{l})d\mu
\\
&\leq C \alpha^{l}R_{i}^{l}\mu(B_{i}).
\end{align*}
Hence $\|b_i\|_l^{l} \leq C\alpha^{l} R_{i}^{l}\mu(B_{i})$.
We invoke the estimate
\begin{equation}\label{eq23}
\left\| \sum_{k\in \mathbb{Z}} \psi(4^k H) \beta_k \right\|_l \lesssim \left\|\left(\sum_{k\in \mathbb{Z}}  |\beta_k|^2\right)^{\frac{1}{2}} \right\|_l.
\end{equation}
Indeed, by duality, this is equivalent to the Littlewood-Paley inequality
$$
\left\| \left(\sum_{k\in \mathbb{Z}} |\psi(4^kH) \beta|^2\right)^{\frac{1}{2}} \right\|_{l'} \lesssim\|\beta\|_{l'}.
$$
This is a consequence of the Gaussian estimates for the kernels of $e^{-tH}$, $t>0$ (this was first proved in \cite{auscher5} using the vector-valued version of the work in \cite{Duong1}.  See \cite{auscher4} or \cite{auscher8} for a more general argument in this spirit or \cite{merdy} for an abstract proof relying on functional calculus).
To apply \eqref{eq23}, observe that the definitions of $r_i$ and $U_i$ yield 
$$
\sum_i U_ib_i = \sum_{k\in \mathbb{Z}} \psi(4^k H) \beta_k
$$
with 
$$\beta_k = \sum_{i, r_i=2^k} \frac{b_i}{r_i}.
$$
Using the bounded overlap property \eqref{eqcsds5}, one has that 
$$
\left\|\left(\sum_{k\in \mathbb{Z}}  |\beta_k|^2\right)^{\frac{1}{2}} \right\|_l^l \leq C\int (\sum_{i}  \frac{|b_i|^l}{r_i^l})d\mu .
$$
Using $R_i\sim r_i$, 
$$
\int (\sum_{i}  \frac{|b_i|^l}{r_i^l})d\mu  \leq  C \alpha^l  \sum_{i}\mu(B_i) .$$
Hence, by  \eqref{eqcsds4}
$$
\mu\left(\left\lbrace x \in M; \bigg|\sum_i U_ib_i(x)\bigg| 
>\frac \alpha 3\right\rbrace\right) \leq  C   \sum_{i}\mu(B_i) \leq \frac{C}{\alpha^l}\int_{M}
(|\nabla f|^l+|V^{\frac{1}{2}}f|^l)d\mu.$$
 Thus, we have obtained
$$\mu\left(\{ x\in M; |H^{\frac{1}{2}}f(x)|>\alpha\}\right)\leq \frac{C}{\alpha^{l}}\int_{M}(|\nabla f|^{l}+V^{\frac{1}{2}}f|^{l})d\mu$$
for all $f\in \Lip(M)\cap \dot{W}_{l,V^\frac{1}{2}}^{1}\cap\dot{W}_{l,V^\frac{1}{2}}^{1}$.
 Moreover, using the density argument of Proposition \ref{DLV} we extend $H^{\frac{1}{2}}$ to a bounded operator acting from $\dot{W}_{l,V^\frac{1}{2}}^{1}$ to $L^{l,\infty}$.
 We already have 
\begin{equation*}
 \| H^{\frac{1}{2}}f\|_{2}\leq  \|\,|\nabla f|\,\|_{2}+\|V^{\frac{1}{2}}f\|_{2}.
\end{equation*}
 Since $V\in A_{\infty}$ implies $V^\frac{1}{2}\in RH_{2}$ --Proposition \ref{CR}--, we see from Corollary \ref{rt8} that
\begin{equation}\label{RRp8}
 \| H^{\frac{1}{2}}f\|_{p}\leq C_{p} \left(\|\,|\nabla f|\,\|_{p}+\|V^{\frac{1}{2}}f\|_{p}\right)
 \end{equation}
 for all $l<p\leq 2$ and $f\in \dot{W}_{p,V}^{1}$.\\
 If $l=1$, we take $1<p<2$. There exists $\epsilon>0$ such that $1<1+\epsilon<p$. The same argument works replacing $l=1$ by $1+\epsilon$. 
 \
 \section{Proof of point 2. of Theorem \ref{RH}}
We first give some estimates for the weak solutions of $-\Delta u+Vu=0$. Then we proceed to a reduction and then give the proof of point 2. of Theorem \ref{RH}.
\subsection{Estimates for weak solutions}\label{sec:weaksol}Let $M$ be a complete Riemannian manifold satisfying $(D)$ and $(P_2)$. Let $B=B(x_{0},R)$ denotes a ball of radius $R>0$ and $u$ a weak solution of $-\Delta u + Vu=0$ in a neighborhood of $\overline{B(x_{0},4R)}$. By a weak solution of $-\Delta u + Vu=0$ in an open set $\Omega$, we mean $u\in L^1_{loc}(\Omega)$ with  $V^{\frac{1}{2}}u,  \nabla u \in L^2_{loc}(\Omega)$ and the equation holds in the distribution sense on $\Omega$. Remark that under the Poincar\'e inequality $(P_{2})$ if $u$ is a weak solution, then $u \in L^2_{loc}(\Omega)$. 
 It should be observed that if $u$ is a weak solution in $\Omega$ of 
$-\Delta u + Vu=0$ then 
\begin{equation}
\label{eq:sub}
 \Delta |u|^2 = 2 V|u|^2 + 2 |\nabla u |^2
\end{equation}
since $\Delta |u|^2 = 2 \langle \Delta u, u\rangle + 2 |\nabla u |^2$
(see \cite{bockner}). 
In particular, $|u|^2$ is  a non-negative subharmonic function in $\Omega$. Hence the lemmas in subsection 3 of section 4 apply to $|u|^2$. In particular
  \begin{equation}
\label{eq:mvi}
\sup_{B(x_{0},R)} |u| \le C(r) \left((|u|^r)_{B(x_{0},\mu R)}\right)^{\frac{1}{r}}
\end{equation}
 holds for any $0<r<\infty$ and $1<\lambda\leq 4$. We have also shown a mean value inequality
against arbitrary $A_{\infty}$ weights.

We state some further estimates that are interesting in their own right assuming 
$V\in A_{\infty}$.  By splitting real and imaginary parts,  we may  suppose $u$ real-valued. All constants are  independent of $B$ and $u$ but they may depend on the constants in the $A_{\infty}$ condition or the $RH_{q}$ condition of $V$ when assumed, on the doubling constant $C_{d}$ and the Poincar\'e inequality $(P_{2})$. 
Let $s$ be any real number such that $\frac{\mu(B)}{\mu(B_{0})}\geq C(\frac{r}{r_{0}})^{s}$ whenever $B=B(x,r),\,x\in B_{0},\,r\leq r_{0}$ ($s=log_{2}C_{d}$ works). 
\\

The proofs of the next 3 lemmas are as in \cite{auscher3}, we skip them. 
\begin{lem}
\label{lem:1}
For all  $1 \leq \lambda < \lambda' \leq 4$ and $k>0$, there is a constant $C$ such that 
$$ 
(|u|^2)_{\lambda B} \leq \frac {C} { (1+ R^{2}V_{B})^k}  (|u|^2)_{\lambda' B}\,.
$$
and
$$ 
 (|\nabla u|^2 + V|u|^2)_{\lambda B} \leq \frac {C} { (1+ R^2V_{B})^k}  (|\nabla u|^2 + V|u|^2)_{\lambda' B}.
$$
\end{lem}

\begin{lem}
\label{lem:2}
For all $1\leq \lambda\leq4$,
 $k>0$,  there is a constant $C$ such that 
$$
(RV_{B})^2  ( |u|^2)_{B} \leq \frac {C} { (1+ R^2V_{B})^k}   ( V|u|^2 )_{\lambda B} .
$$ 
\end{lem}

\begin{lem}
\label{lem:3} 
For all 
$1<\lambda\leq 4$,
 $k>0$  and $\max(s,2) <r<\infty$, there is a constant $C$ such that 
$$
(RV_{B})^2  (|u|^2)_{B} \leq  \frac {C} { (1+ R^2V_{B})^k}    ( |\nabla u|^r)_{\lambda B}^{\frac{2}{r}}.
$$
\end{lem}
The main tools to prove these lemmas are the improved Fefferman-Phong inequality of Lemma \ref{lem:FP}, the Caccioppoli type inequality which holds on complete Riemannian manifolds, Poincar\'{e} inequality, subharmonicity of $|u|^{2}$, Lemma \ref{subw} and the Morrey embedding theorem with exponent $\alpha=1-\frac{s}{r}$ (\cite{hajlasz4}, Theorem 5.1, p. 23) to prove Lemma \ref{lem:3}.

For the remaining lemmas, we moreover assume that $M$ is of polynomial type: every ball $B$ of radius $r>0$ satisfies \begin{equation*}\tag{$L_{\sigma}$}
\mu(B)\geq cr^{\sigma},
\end{equation*} 
 and
\begin{equation*}\tag{$U_{\sigma}$}
\mu(B)\leq Cr^{\sigma} 
\end{equation*}
with $\sigma=d$ if $r\leq 1$ and $\sigma=D$ for $r\geq  1$ and $d\leq D$.
Note that if $(L_{\sigma})$ holds then $\sigma \geq n$ where $n$ is the topological dimension of $M$ (see \cite{saloff3}). Recall that under $(L_{\sigma})$ and $(U_{\sigma})$, $s=D$ works and that $\mu(B(x,r))\geq cr^{\lambda}$ for all $r>0$ with any $\lambda \in [d,D]$. We also recall that the exponent $p_0$ is that appearing in Proposition \ref{RTN}.

\begin{lem}
\label{lem:4}  
Assume $V \in RH_{q}$. Let $B$ be a ball of radius $R>0$ and $\sigma=d$ if $R\leq 1$ and $\sigma=D$ if $R\geq 1$. Set $\tilde q = q^{*}_{\sigma}$ if $q^{*}_{\sigma}<p_{0}$ ($\frac{1}{q^{*}_{\sigma}}=\frac{1}{q}-\frac{1}{\sigma}$) and $\tilde{q}$ arbitrary in $]2,\,p_{0}[$ if not. Then
 for all $k>0$ there is a constant $C=C(\sigma)$ independent of $B$ such that
$$
   \left( (|\nabla u|^{\tilde q})_{B} \right)^{\frac{1}{\tilde q}} \leq   \frac {C} { (1+ R^2V_{B})^k}  \left(( |\nabla u|^2 + V|u|^2 )_{4 B}\right)^{\frac{1}{2}}. 
$$
\end{lem}
\begin{lem}
\label{lem:5}
Assume $V \in RH_{q}$ with $\frac{D}{2}\leq q<\frac{p_0}{2}$. Let $B$ be a ball of radius $R>0$ and $\sigma=d$ if $R\leq 1$ and $\sigma=D$ if $R\geq 1$. Set $\tilde{q}=q^{*}_{\sigma}$ if $q^{*}_{\sigma}<p_{0}$ and $\tilde{q}$ arbitrary in $]2q,\,p_{0}[$ if not. 
 Then, there is a constant $C=C(\sigma)$ such that
$$
\left(( |\nabla u|^{ \tilde{q}})_{B} \right)^{\frac{1}{ \tilde{q}}} \leq  C \left(( |\nabla u|^2)_{4 B}\right)^{\frac{1}{2}}, 
$$
\end{lem}
 We give the proofs of Lemma \ref{lem:4} and \ref{lem:5} since they are not exactly the same as the one in the Euclidean case. Before the proof of Lemma \ref{lem:4}, we need the following theorem for the boundedness of the Riesz potential.
\begin{thm}(\cite{chen})\label{RP8} Let $M$ be a complete Riemannian manifold satisfying $(D)$ and $(P_{2})$. Moreover, assume that $M$ satisfies 
\begin{equation*}\tag{$L_{\lambda}$}
\mu(B)\geq cr^{\lambda}
\end{equation*} 
for every $x\in M$ and $r>0$.
Then $(-\Delta)^{-\frac{1}{2}}$ is $L^{p}-L^{p*}$ bounded with $1<p,\,p^{*}<\infty$ and $p^{*}=\frac{\lambda p}{\lambda-p}$, that is,
$$ 
\|(-\Delta)^{-\frac{1}{2}}f\|_{p^{*}}\leq C(p,\lambda)\|f\|_{p}.
$$ 
\end{thm}
\begin{proof} In \cite{chen}, Chen proves this theorem for Riemannian manifolds with non-negative Ricci curvature. His proof still works under our hypotheses. The properties that he used for these manifolds are first the lower and upper gaussian estimates for the heat kernel which holds on Riemmanian manifolds satisfying $(D)$ and $(P_2)$. Secondly, he applied an argument from the proof of the $L^{p}-L^{p^{*}}$ boundedness of the Riesz potential in the Euclidean case (\cite{stein2}, Chapter V, Theorem 1) which remains true since we have $(D)$, $(P_{2})$ and $(L_{\lambda})$ with $\lambda\geq n=dim\, M$. 
\end{proof}
\begin{proof}[Proof of Lemma \ref{lem:4}]
First note that if $q\le \frac{2\sigma}{\sigma+2}$ then $\tilde q\le 2$ and the conclusion (useless for us) follows by a mere H\"older inequality. Henceforth, we assume $q>\frac{2\sigma}{\sigma+2}$. Also, by Lemma \ref{lem:1}, it suffices to obtain the estimate with $k=0$. Let us assume $\mu=4$ for simplicity of the argument.
Let $v$ be the harmonic function  on $ 4 B$ with $v=u$ on $\partial( 4 B)$ and set $w=u-v$ on $ 4 B$. Since $w=0$ on $\partial( 4B)$, the fact that an harmonic function minimises Dirichlet integral among functions with the same boundary implies 
$$
 (\aver{ 4 B}   |\nabla w|^2 \big)^{\frac{1}{2}}  \le 
2(\aver{ 4 B}   |\nabla u|^2 \big)^{\frac{1}{2}} .
$$
 By the elliptic estimate for the harmonic function $v$ (\cite{auscher1}, Theorem 2.1), we have for $p<p_0$  
  \begin{equation}\label{VU}
\big ( \aver{B} |\nabla v|^p \big)^{\frac{1}{p}} \le C (\aver{ 4 B}   |\nabla v|^2 \big)^{\frac{1}{2}}  \le 2 C 
(\aver{ 4 B}   |\nabla u|^2 \big)^{\frac{1}{2}}.
\end{equation}
 Let $1<\nu<\lambda <4$ and 
 $\eta$ be a smooth non-negative function, bounded by 1, equal to 1 on $\nu B$ with support contained in $\lambda B$ and whose gradient is bounded by $\frac{C}{R}$.  As $\Delta w= \Delta u = V u $ on $4B$, we have
$$
\Delta(w\eta) =  Vu \eta + \nabla w\cdot \nabla \eta  +div( w \nabla \eta)   \quad {\rm on \ } M.
$$
It comes that
\begin{align*}
\nabla (w\eta)(x) &=\nabla(-\Delta)^{-1}(-\Delta)(w\eta)(x)
\\&=
\nabla(-\Delta)^{-\frac{1}{2}}(-\Delta)^{-\frac{1}{2}}(-Vu\eta)(x)+\nabla(-\Delta)^{-\frac{1}{2}}(-\Delta)^{-\frac{1}{2}}(-\nabla w.\nabla \eta)(x)
\\
&+\nabla(-\Delta)^{-1}(-div(w\nabla \eta))(x)
\\
&= I+II+III.
\end{align*}
Let us begin with $$
III=\nabla(-\Delta)^{-\frac{1}{2}}(-\Delta)^{-\frac{1}{2}}div(-w\nabla \eta)(x)=(\nabla (-\Delta)^{-\frac{1}{2}})(\nabla(-\Delta)^{-\frac{1}{2}})^{*}(-w\nabla \eta)(x).
$$
 Let $\eta'$ be a smooth function, bounded by $1$, equal to $1$ on $\lambda B$ with support contained in $\lambda'B$ with $\lambda'<4$ and whose gradient is bounded by $\frac{C}{R}$.
The Riesz transform $\nabla (-\Delta)^{-\frac{1}{2}}$ is $L^{p}(M)$ bounded for $1<p<p_{0}$. By duality, $(\nabla(-\Delta)^{-\frac{1}{2}})^{*}$ is $L^{p}(M)$ bounded for $p'_0<p<\infty$. Hence for $2<p<p_0$
\begin{align*}
\left(\int_{M}|III|^{p}d\mu\right)^{\frac{1}{p}}&\leq C\left(\int_{M}|w\eta '|^{p}|\nabla \eta|^{p}d\mu\right)^{\frac{1}{p}}
\\
&\leq \frac{C}{R}\left(\int_{M}|\nabla (w\eta')|^{p_{*}}d\mu\right)^{\frac{1}{p_{*}}}.
\end{align*}
We used the Sobolev inequality which holds under $(D)$, $(P_{2})$ and $\mu(B(x,r))\geq cr^{\sigma}$ for all $r>0$ with $p_{*\sigma}<p$ defined by $p_{*\sigma}=\frac{\sigma p}{\sigma +p}$ that is $(p_{*})^{*}=p$ (see \cite{saloff3}).
 We use the $L^{q}-\,L^{q^{*}_{\sigma}}$ boundedness of the Riesz potential $(-\Delta)^{-\frac{1}{2}}$ and the  $L^{p}$ boundedness of the Riesz transform $\nabla (-\Delta)^{-\frac{1}{2}}$ for $1<p<p_{0}$ to get the estimates for II and I. First for II, we have for all $2\leq p<p_0$
\begin{align*}
\left(\int_{M}|II|^{p}d\mu\right)^{\frac{1}{p}}&\leq C\left(\int_{M}|(-\Delta)^{-\frac{1}{2}}(\nabla (w\eta').\nabla \eta)|^{p}d\mu\right)^{\frac{1}{p}}
\\
&\leq \frac{C}{R}\left(\int_{M}|\nabla (w\eta')|^{p_{*\sigma}}d\mu\right)^{\frac{1}{p_{*\sigma}}}
 \\
 &\leq \frac{C}{R}\left(\int_{M}|\nabla (w\eta')|^{p_{*\sigma}}d\mu\right)^{\frac{1}{p_{*\sigma}}}
 \\
 &=\frac{C}{R}\left(\int_{M}|\nabla (w\eta')|^{p_{*\sigma}}d\mu\right)^{\frac{1}{p_{*\sigma}}}.
 \end{align*}
 Now, it remains to look at I. Take $p=q^{*}_{\sigma}$ if $q^{*}_{\sigma}<p_{0}$ and if not any $2<p<p_{0}$. It follows that
\begin{align*}
\left(\int_{M}|I|^{p}d\mu\right)^{\frac{1}{p}}&\leq C\bigg( \int_{
M} |Vu\eta|^{p_{*\sigma}} d\mu\bigg)^{\frac{1}{p_{*\sigma}}}  \le C\mu(B)^{\frac{1}{p_{*\sigma}}}\bigg( \aver{
\lambda B} |V|^qd\mu \bigg)^{\frac{1}{q}} \sup_{\mu B } |u|
\end{align*}
since $p_{*\sigma}\leq q$ in the two cases.
Using the $RH_{q}$ condition on $V$, we obtain
\begin{equation}
\label{eq:I}
\big (\int_{M} |I|^{p}d\mu\big)^{\frac{1}{p}} \le C\mu(B)^{\frac{1}{p_{*\sigma}}}  \, \aver{\lambda B} Vd\mu \sup_{\mu B } |u|.
\end{equation}
Now, if $\lambda<\gamma<4$, the subharmonicity of $|u|^2$ and Lemma \ref{subw} yield
$$
 \aver{\lambda B} V d\mu\,\sup_{\lambda B}|u| \le  C   \aver{\gamma B} V d\mu \,  \big(\aver{\gamma B} |u|^2d\mu \big)^{\frac{1}{2}}.
$$
 It follows from Lemma \ref{lem:2} and $(U_{\sigma})$ that 
 $
\big (\int_{M} I^{p}d\mu\big)^{\frac{1}{p}}\leq C\mu(B)^{\frac{1}{p}}  \big(  \avert{4 B}  V|u|^2d\mu\big)^{\frac{1}{2}}.
$
Therefore, we showed that
$$
\left(\int_{M}|\nabla (w\eta)|^{p}d\mu\right)^{\frac{1}{p}}\leq \frac{C}{R}\left(\int_{M}|\nabla (w\eta')|^{p_{*}}d\mu\right)^{\frac{1}{p_{*}}}+
C\mu(B)^{\frac{1}{p}} \left(  \aver{4 B} V|u|^2d\mu\right)^{\frac{1}{2}}.
$$
We repeat the same process 
and after a finite iteration ($K=(\sigma [\frac{1}{2}-\frac{1}{p}]+1$) times), using $(U_{\sigma})$ we get 
$$\left(\aver{B}|\nabla w|^{\tilde{q}}d\mu\right)^{\frac{1}{\tilde{q}}}\leq C\left(\aver{4 B}|\nabla w|^{2}d\mu\right)^{2}+C \left(  \aver{4 B}  V|u|^2 d\mu \right)^{\frac{1}{2}}.
$$
We derive therefore the desired inequality for $\nabla u$ from the estimates obtained for $\nabla v$ and $\nabla w$.
\end{proof}

\begin{proof}[Proof of Lemma \ref{lem:5}] Since $V\in RH_q$ and $q\geq \frac{D}{2}$, we may assume $q>\frac{D}{2}$ by self-improvement. Let $\sigma=d$ if $R\leq 1$ and $\sigma=D$ if $R\geq1$. We apply the same arguments as in the proof of the previous lemma. The only difference is that since $2q>s=D$, we use Lemma \ref{lem:3} with $k=0$, $r=2q$, and $s=D$ instead of Lemma \ref{lem:2} in the estimate for the term I. We then obtain 
\begin{equation}\label{gradtild}
\big( \aver{B} |\nabla u|^{\tilde{q}} \big)^{\frac{1}{ \tilde{q} }} \le C \big( \aver{4B} |\nabla u|^{ 2q} \big)^{\frac{1}{ 2q }}
\end{equation}
where $p=q^{*}_{\sigma}$ if $q^{*}_{\sigma}<p_0$ and if not we take any $2<p<p_0$.
Since $2q<p_0$, if we take $p=\tilde{q}\in ]2q,p_0[$ in (\ref{gradtild}) we can apply Lemma \ref{buck2}  and improve the exponent $2q$ to $2$. Thus, we get 
$$
\big( \aver{B} |\nabla u|^{ {\tilde{q}}} \big)^{\frac{1}{ \tilde{q} }} \le C \big( \aver{4 B} |\nabla u|^{ 2} \big)^{\frac{1}{ 2 }}
$$
Remark that when $q>D$, $q^{*}_{\sigma}=\infty$ and therefore we have our lemma for any $2q< p<p_{0}$.
\end{proof}

\subsection{A reduction}

It is sufficient to prove the  $L^p$ boundedness of $\nabla H^{-\frac{1}{2}}$ and of $V^{\frac{1}{2}}H^{-\frac{1}{2}}$ for the appropriate range of $p$. As we have seen in the introduction, the case $1<p\leq 2$ does not need any assumption on $V$. We henceforth assume $p>2$ and $V\in A_{\infty}$.

By duality, we know that $H^{-\frac{1}{2}}div$ and $H^{-\frac{1}{2}}V^{\frac{1}{2}}$ are bounded on $L^p$ for $2<p<\infty$. Thus, if $\nabla H^{-\frac{1}{2}}$ is also bounded on $L^p$. It follows that $\nabla H^{-1} div$ and $\nabla H^{-1}V^{\frac{1}{2}}$ are bounded on $L^p$. 

Reciprocally, if $\nabla H^{-1} div$ and $\nabla H^{-1}V^{\frac{1}{2}}$ are bounded on $L^p$, then their adjoints are bounded on $L^{p'}$. Thus, if $F\in C_{0}^\infty(M, TM)$, 
\begin{align*}
\|H^{-\frac{1}{2}} div F\|_{p'} &= \|H^{\frac{1}{2}} H^{-1}div F\|_{p'}\\
&\leq C (\|\,|\nabla H^{-1} div F|\,\|_{p'}
+ \| V^{\frac{1}{2}} H^{-1} div F\|_{p'}) \leq C \|F\|_{p'}
\end{align*}
where the first inequality follows from Theorem \ref{RRT}. By duality, we have that $\nabla H^{-\frac{1}{2}}$ is bounded on $L^p$. 

The same treatment can be done on $V^{\frac{1}{2}} H^{-\frac{1}{2}}$.  We have obtained 
\begin{lem} Let $M$ be a complete Riemannian manifold. If $V\in A_{\infty}$ and $p>2$,  the $L^p$ boundedness of  $\nabla H^{-\frac{1}{2}}$  is equivalent to that of $\nabla H^{-1} div$ and $\nabla H^{-1}V^{\frac{1}{2}}$, and  the $L^p$ boundedness of  $V^{\frac{1}{2}} H^{-\frac{1}{2}}$  is equivalent to that  of  $V^{\frac{1}{2}} H^{-1} V^{\frac{1}{2}}$ and $V^{\frac{1}{2}} H^{-1} div$.
\end{lem}
 Hence, to prove point 2. of Theorem \ref{RH}, it suffices the $L^{p}$ boundedness of the operators: $\nabla H^{-1} div,\, \nabla H^{-1}V^{\frac{1}{2}}$, $V^{\frac{1}{2}} H^{-1}V^{\frac{1}{2}}, V^{\frac{1}{2}} H^{-1} div$.
\subsection{Proof of point 2. of Theorem \ref{RH}}
\begin{prop} \label{thC'}Let $M$ be a complete Riemannian manifold satisfying $(D)$ and $(P_2)$. Assume that $V \in RH_q$ for some $q>1$. Then for 
$2<p< 2(q+\epsilon)$, for some $\epsilon>0$ depending only on $V$,  $f\in C^\infty_{0}(M,\mathbb{C})$ and $F\in C_{0}^\infty(M,TM)$, 
$$
\|V^{\frac{1}{2}} H^{-1} V^{\frac{1}{2}}f\|_p  \leq C_{p} \|f\|_{p}, \quad  \|V^{\frac{1}{2}} H^{-1} div F\|_p \leq C_p\|F\|_p.
$$
\end{prop}

\begin{prop} \label{thD'} Let $M$ be a complete Riemannian manifold of polynomial type satisfying $(P_2)$. 
Let $V \in RH_q$  for some $q>1$. 
If $q^{*}_{D}<p_0$, let $p=q^{*}_{D}$. If $q^{*}_{D}\geq p_0$, we take any $2<p<p_{0}$. Then for all $f\in C^\infty_{0}(M,C)$ and $F\in C_{0}^\infty(M, TM)$, 
$$
\|\nabla  H^{-1} V^{\frac{1}{2}}f\|_p  \leq C_{p} \|f\|_{p}, \quad  \|\,|\nabla H^{-1} div F|\,\|_p \leq C_p\|F\|_p.
$$
\end{prop}
 
 The interest of such a reduction is that this allows us to use properties of weak solutions of  $H$.
 Note that Proposition \ref{thD'} is void if $q\leq \frac {2D}{D+2} $ as $q^*_{D}\leq 2$. Note also that $q^*_{D}<2q$ exactly when $q<\frac{D}{2}$. In this case, this statement yields a smaller range than 
 the interpolation method in Section 6.
\begin{proof}[Proof of Proposition \ref{thC'}]
 Fix a ball $B=B(x_{0},R)$ and let $f\in C^\infty_{0}(M)$ supported away from $\overline{4B}$. Then $u=H^{-1}V^{\frac{1}{2}}f$ is well defined on $M$ with 
$\|V^{\frac{1}{2}} u\|_{2}+ \|\,|\nabla u|\, \|_{2} \leq \|f\|_{2}$ by construction of $H$ and 
$$ 
\int_{M} (V u \varphi + \nabla u \cdot \nabla \varphi)d\mu = \int_{M} V^{\frac{1}{2}} f \varphi d\mu
$$ for all $\varphi\in L^2(M)$ with $\|V^{\frac{1}{2}} \varphi\|_{2}+ \|\,| \nabla\varphi|\, \|_{2}<\infty$.
In particular, the support condition on $f$ implies that $u$ is a weak solution of $-\Delta u + Vu=0$ in a neighborhood of $\overline{4B}$, hence $|u|^2$ is subharmonic there. Let $r$ such that $V\in RH_{r}$. Note that by Proposition \ref{CR},  $V^{\frac{1}{2}} \in RH_{2r}$. From Corollary \ref{lem:rh2q} with $V^{\frac{1}{2}}$, $|u|^2$ and $s=\frac{1}{2}$, we get
$$
\big(\aver{B} (V^{\frac{1}{2}}|u| )^{2r}d\mu \big)^{\frac{1}{2r}} \leq C \, \aver{4B} V^{\frac{1}{2}} |u| d\mu.
$$
Thus, \eqref{T:shen} holds with $T=V^{\frac{1}{2}}H^{-1} V^{\frac{1}{2}}$, $q_{0}=2r$, $p_{0}=2$ and $S=0$. By Theorem \ref{theor:shen}, 
$V^{\frac{1}{2}}H^{-1} V^{\frac{1}{2}}$ is bounded on $L^p$ for $2<p<2r$. 

The argument is the same  for $V^{\frac{1}{2}}H^{-1}div$. This finishes the proof. 
\end{proof}

\begin{proof}[Proof of Proposition \ref{thD'}]
We assume $q>\frac{2D}{D+2}$, that is $q^{*}_{D}>2$, otherwise there is nothing to prove. We consider first the operator $\nabla H^{-1} V^{\frac{1}{2}}$.

 Assume $q<\frac{D}{2}$. 
Fix a ball $B$  of radius $R$ and let $f\in C^\infty_{0}(M)$  supported away from $\overline{4B}$. Let $u=H^{-1}V^{\frac{1}{2}}f$. As before,  the support condition on $f$ implies that $u$ is a  weak solution of $-\Delta u + Vu=0$ in a neighborhood of $\overline{4B}$.  Thanks to Lemma \ref{lem:4}, 
\eqref{T:shen} holds with $T=\nabla H^{-1} V^{\frac{1}{2}}$, $q_{0}= q^{*}_{D}\leq q^{*}_{d}$ if $q^{*}_{D}<p_{0}$ and if not $q_{0}=p_{0}-\epsilon'$ for any $\epsilon'>0$, and $S=\left(\mathcal{M}(|V^{\frac{1}{2}}H^{-1}V^{\frac{1}{2}}|^{2})\right)^{\frac{1}{2}}$.
The maximal theorem --Theorem \ref{MIT}-- and Proposition \ref{thC'} show that $S$ is bounded on $L^{p}(M)$ for $1<p<2q$. Then Theorem \ref{theor:shen} implies that 
$\nabla H^{-1} V^{\frac{1}{2}}$ is bounded on $L^p(M)$ for $2<p < p_{0}$ if $q^{*}_{D}\geq p_{0}$. If $q^{*}_{D}<p_{0}$, by the self-improvement of reverse H\"older estimates we can replace $q$ by a slightly larger value and, therefore we get the $L^p$ boundedness of $\nabla H^{-1}V^{\frac{1}{2}}$  for $p\leq q^{*}_{D}$. 

Assume next that  $\frac{D}{2}\leq q<D$ and $2q<p_0$. Again, we may as well assume $q>\frac{D}{2}$. In this case $q^{*}_{D}>2q$. Then, Lemma \ref{lem:5}  yields, this time, \eqref{T:shen}  with $T=\nabla H^{-1} V^{\frac{1}{2}}$, $q_{0}=q^{*}_{D}$ if $q^{*}_{D}<p_0$ and if not $q_0=p_{0}-\epsilon'$ for any $0<\epsilon'<p_0-2q$, and $S=0$. Theorem \ref{theor:shen} asserts that 
$\nabla H^{-1} V^{\frac{1}{2}}$ is bounded on $L^p$ for $2<p< p_{0}$ if $q^{*}_{D}\geq p_{0}$ and, by the self-improvement
of the $RH_{q}$ condition, it holds for $p\leq q^{*}_{D}$ if $q^{*}_{D}<p_{0}$.

Finally, if $q\geq D $, then Lemma \ref{lem:5} yields \eqref{T:shen} for any $2<q_{0}<p_{0}$  with $T=\nabla H^{-1} V^{\frac{1}{2}}$ and $S=0$. Theorem \ref{theor:shen} shows then that
$\nabla H^{-1} V^{\frac{1}{2}}$ is bounded on $L^p$ for $2<p < p_{0}$. 

The argument is the same  for $\nabla H^{-1}div$ and the proof is therefore complete.
\end{proof}

\section{Case of Lie groups} Consider $G$ a  simply connected Lie group. Assume that $G$ is unimodular and let $d\mu$ be a fixed Haar measure on $G$. Let $X_{1},...,X_{k}$ be a family of left invariant vector fields such that the $X_{i}$'s satisfy a H\"{o}rmander condition. In this case the Carnot-Carath\'{e}odory metric $\rho$ is a distance, and the metric space $(G,\rho)$ is complete and has the same topology as $G$ as a manifold (see \cite{coulhon6} page 1148). Denote $V(r)=\mu(B(x,r))$ for all $x\in G$. An important result of Guivarc'h \cite{guivarch} says that, either there exists an integer $D$ such that $cr^{D}\leq V(r)\leq Cr^{D}$  for all $r>1$, or $e^{cr}\leq V(r)\leq Ce^{Cr}$ for all $r>1$ with $V(r)=\mu(B(x,r))=\mu(B(y,r))$, for all $x,\,y\in G$ and $r>0$.
In the first case we say that $G$ has polynomial growth, while in the second case $G$ has exponential growth.
For small $r$, a result of \cite{nagel} implies that there exists an integer $d$ such that $cr^{d}\leq V(r)\leq Cr^{d}$ for $0<r<1$. 
Suppose that $G$ has polynomial growth.
Then there exists $C_{1}>0$ such that
\begin{equation}\label{LG}
C_{1}^{-1}r^{d}\leq V(r)\leq C_1r^{d},\quad  0\leq r\leq 1,
\end{equation}
\begin{equation}\label{UG}
C_{1}^{-1}r^{D}\leq V(r)\leq C_1r^{D},\quad  1\leq r<\infty.
\end{equation}
We say that $d$ is the local dimension of $G$ and $D$ is the dimension at infinity. We assume that $d\geq 3$ and $d\leq D$ --If $G$ is nilpotent and since $G$ is simply connected, we have $d\leq D$ (see \cite{coulhon7})--. 
In particular $(D)$ holds with $s=D$.
Moreover $G$ satisfies a Poincar\'{e} inequality $(P_{1})$: there exists $C>0$ such that for all ball $B$ of radius $r>0$ we have for every smooth function $u$,
 \begin{equation}\tag{$P_{1}$}
 \int_{B}|u-u_{B}|d\mu \leq Cr\int_{2B}|Xu|d\mu 
 \end{equation}
 (see \cite{saloff2}, \cite{varopoulos2}) where $|Xu|=\left(\sum_{i=1}^{k}|X_{i}u|^{2}\right)^{\frac{1}{2}}$.

For the rest of this section, we consider $G$ a Lie group as above with polynomial growth and set $\Delta=\sum_{i=1}^kX_{i}^2$.

Let us check the validity of our approach to obtain Theorem \ref{th:1}, Theorem \ref{RH} and Theorem \ref{RRT} for $G$. The main tools used to prove those theorems still hold:
\begin{itemize}
\item The Riesz transform $\nabla (-\Delta)^{-\frac{1}{2}}$ is $L^p$ bounded for all $1<p<\infty$. This result was proved by Alexopoulos \cite{alexopoulos}.

\item An improved Fefferman-Phong inequality of type (\ref{eq:FP}) holds on $G$ with $\beta=\frac{p}{p+D(\alpha-1)}$.
\item We get a Calder\'{o}n-Zygmund decomposition analogous to that of Proposition \ref{CZ}. Thanks to this decomposition, we get the analog of Theorem \ref{RRT} as in section 7.
\item Theorem \ref{HI} proved in section 6 remains true for Lie groups with polynomial growth (we use the same proof). 
\item The argument of complex interpolation (valid on $G$) allows us to obtain Theorem \ref{RH} part 1.

\item Let $u$ a weak solution of $-\Delta u+Vu=0$ on $G$, then $u$ satisfies some mean values inequalities as in Lemma \ref{Eu}, \ref{subw} and Corollary \ref{lem:rh2q}. We mention that the analogous of Lemma \ref{Eu} was proved by Li \cite{li}, \cite{li3} for nilpotents groups using estimations for the heat kernel and its first and second derivatives. 
\item The lemmas in section \ref{sec:weaksol} still hold in our case: $G$ is of polynomial type. The Sobolev inequality and the Morrey embedding --with $\alpha=1-\frac{n}{p}$ and $1-\frac{n}{p}\notin \mathbb{N}$-- hold for any $n\in [d,D]$  (see Theorem VIII.2.10 of \cite{coulhon7}).
We also have that $\Delta^{-\frac{1}{2}}$ is bounded from $L^{p}$ to $L^{\frac{np}{n-p}}$ for any $n\in[d,D]$ and $p<n$ (Theorem VIII.2.3 of \cite{coulhon7}). Thus we get similar lemmas to that of section \ref{sec:weaksol} this time on a Lie group $G$ of polynomial growth. 
\end{itemize} With all these ingredients, we establish the following theorem analog to Theorem \ref{RH}.

 \begin{thm} \label{RHG} Let $G$ be a simply connected Lie group with polynomial growth and assume $3\leq d\leq D$. Let $V\in RH_{q}$ for some $q>1$.
\begin{itemize} 
\item[1.] Then  for any smooth function $u$,
\begin{equation}\label{RpG}
\|\,|\nabla
u|\,\|_{p}+\|V^{\frac{1}{2}}u\|_{p}\lesssim\|(-\Delta+V)^{\frac{1}{2}}u\|_{p}\;\textrm{
for } 1<p<2(q+\epsilon).
\end{equation}
\item[2.] Assume $q\geq\frac{D}{2}$. Consider
   \begin{equation}\label{gradg}
 \|\,|\nabla u|\,\|_{p}\lesssim \|(-\Delta+V)^{\frac{1}{2}}u\|_{p}
 \end{equation} 
 for all smooth function $u$.
\begin{itemize}
\item[a.] if $\frac{D}{2}<q<D$, (\ref{gradg}) holds for $1<p<q^{*}_{D}+\epsilon$,
\item[b.] if  $q\geq D$, (\ref{gradg}) holds for $1<p<\infty$ .
\\
\end{itemize}
\end{itemize}
 \end{thm}

\begin{rem} Li \cite{li}, \cite{li3} proved point 2. of Theorem \ref{RHG} if $G$ is in addition Nilpotent. 
\end{rem}

\bibliographystyle{plain}
\bibliography{latex2}
\end{document}